\documentclass[11pt,oneside,english]{amsart}
\usepackage[latin9]{inputenc}
\usepackage{babel}
\usepackage{float}
\usepackage{bm}
\usepackage{amsthm}
\usepackage{amssymb}
\usepackage{stmaryrd}
\usepackage{graphicx}
\usepackage[unicode=true,pdfusetitle,
 bookmarks=true,bookmarksnumbered=false,bookmarksopen=false,
 breaklinks=false,pdfborder={0 0 1},backref=false,colorlinks=false]
 {hyperref}

\makeatletter

\providecommand{\tabularnewline}{\\}

\numberwithin{equation}{section}
\numberwithin{figure}{section}
\theoremstyle{plain}
\newtheorem{thm}{\protect\theoremname}[section]
\theoremstyle{definition}
\newtheorem{defn}[thm]{\protect\definitionname}
\theoremstyle{remark}
\newtheorem{rem}[thm]{\protect\remarkname}
\theoremstyle{plain}
\newtheorem{lem}[thm]{\protect\lemmaname}
\theoremstyle{plain}
\newtheorem{prop}[thm]{\protect\propositionname}
\theoremstyle{definition}
\newtheorem{example}[thm]{\protect\examplename}

\usepackage{amsfonts}
\usepackage{amsthm}
\usepackage{cite}
\usepackage{url}
\usepackage{pb-diagram}

\address{University of Arizona, Department of Mathematics, 617 N Santa Rita
Ave, Tucson, AZ 85721, United States.}
\email{ajzerouali@math.arizona.edu}

\makeatother

\providecommand{\definitionname}{Definition}
\providecommand{\examplename}{Example}
\providecommand{\lemmaname}{Lemma}
\providecommand{\propositionname}{Proposition}
\providecommand{\remarkname}{Remark}
\providecommand{\theoremname}{Theorem}

\begin{document}
\title{Twisted Conjugation on Connected Simple Lie Groups and Twining Characters}
\author{Ahmed J. Zerouali}
\begin{abstract}
This article discusses the twisted adjoint action $\mathrm{Ad}_{g}^{\kappa}:G\rightarrow G$,
$x\mapsto gx\kappa(g^{-1})$ given by a Dynkin diagram automorphism
$\kappa\in\mathrm{Aut}(G)$, where $G$ is compact, connected, simply
connected and simple. The first aim is to recover the classification
of $\kappa$-twisted conjugacy classes by elementary means, without
invoking the non-connected group $G\rtimes\langle\kappa\rangle$.
The second objective is to highlight several properties of the so-called
\textit{twining characters} $\tilde{\chi}^{(\kappa)}:G\rightarrow\mathbb{C}$,
as defined by Fuchs, Schellekens and Schweigert. These class functions
generalize the usual characters, and define $\kappa$-twisted versions
$\tilde{R}^{(\kappa)}(G)$ and $\tilde{R}_{k}^{(\kappa)}(G)$ ($k\in\mathbb{Z}_{>0}$)
of the representation and fusion rings associated to $G$. In particular,
the latter are shown to be isomorphic to the representation and fusion
rings of the \textit{orbit Lie group} $G_{(\kappa)}$, a simply connected
group obtained from $\kappa$ and the root data of $G$.
\end{abstract}

\maketitle

\section{\textbf{Introduction}}

Let $G$ be a compact connected semi-simple Lie group, and $\kappa\in\mathrm{Aut}(G)$
an element induced by a Dynkin diagram automorphism of $G$. The $\kappa$-\textit{twisted
adjoint action} of $G$ on itself is defined as:
\[
\mathrm{Ad}_{g}^{\kappa}:G\longrightarrow G\mbox{, }x\longmapsto gx\kappa(g^{-1})\mbox{, }g\in G.
\]
Early publications discussing this action are Gantmacher's work on
automorphisms of complex groups \cite{[Ga39]}, and de Siebenthal's
study of compact non-connected Lie groups \cite{[Sib56]}. Since then,
twisted conjugation appeared sporadically in the literature, mostly
in the context of representations of non-connected groups and of twisted
affine Lie algebras \cite{[Seg68], [Stein68], [Jant73], [CaSGLT72]}.
Many results in this direction can be found in Mohrdieck's work on
twisted conjugacy classes \cite{[Mo00],[Mo03]}, in Wendt's study
of orbital integrals \cite{[Wen01]}, as well as in Springer's note
\cite{[Spr06]}.

Aside from representation theory, twisted conjugation has found several
applications in mathematical physics and gauge theory in the last
two decades. One particularly relevant paper is Fuchs, Schellekens
and Schweigert's \cite{[FSS96]} on \textit{twining characters}, which
introduced twisted analogs of the usual characters and recovered an
important character formula of Jantzen, in order to solve the so-called
``resolution of fixed points'' problem in conformal field theory.
Later, Mohrdieck and Wendt studied twisted conjugacy classes as D-branes
for the Wess-Zumino-Witten model in \cite{[MW04]} (see also \cite{[CW08]}),
as well as in the context of moduli spaces of principal bundles in
\cite{[MW09]}. More recently, Hong addressed these questions in \cite{[Hon17a],[Hon17b], [HK18]},
where he constructed an analog of the fusion ring for twining characters
and developed a notion of twisted conformal blocks with Kumar.

Our interest in twisted conjugation stems from the study of certain
symplectic invariants of \textit{twisted quasi-Hamiltonian spaces},
as discussed in work of Boalch and Yamakawa \cite{[BY15]}, Knop \cite{[Kno16]}
and Meinrenken \cite{[Mein17]}. These are finite-dimensional models
for symplectic Banach manifolds endowed with Hamiltonian actions of
twisted loop groups. We study these spaces in \cite{[Zer18b]}, where
we use some of the results of the present work.

In light of these developments, the purpose of this paper is twofold:
First, review twisted conjugation with direct proofs of the main structure
results, and second, highlight certain properties of twining characters
from this perspective. Some features of our approach are that: \begin{itemize}

\item[(i)] The action $\mathrm{Ad}^{\kappa}$ is studied as an intrinsic
action of $G$ on itself, without resorting to the theory of non-connected
groups.

\item[(ii)] The classification of twisted conjugacy classes is recovered
by elementary means, from a differential-geometric perspective.

\item[(iii)] Twining characters are viewed as functions on the group
$G$, as opposed to formal characters of an affine Kac-Moody algebra.

\item[(iv)] The level $k$ fusion ring for twining characters is
constructed without invoking twisted loop groups or affine Lie algebras.

\end{itemize}

Keeping the notation above, suppose henceforth that $G$ is simply
connected and simple. Let $T\subseteq G$ be a maximal torus with
fixed-point subgroup $T^{\kappa}\subseteq T$, and denote by $\mathfrak{R}=\mathfrak{R}(G,T)$
the root system of $G$. The upcoming sections are organized as follows.
In section \ref{Sec_Diag_Aut_Root_Systs}, we introduce the notations
and associate two root systems to the data $(G,\kappa)$: the \textit{folded
root system} $p(\mathfrak{R})$ and the \textit{orbit root system}
$\mathfrak{R}_{(\kappa)}$, both supported in the dual of $\mathfrak{t}^{\kappa}=\mathrm{Lie}(T^{\kappa})$.
We then gather the properties of these root systems that form the
basis of sections \ref{Sec_Twisted_Conj_Classes} and \ref{Sec_Twining_Char}.

The structure of $\kappa$-twisted conjugacy classes in $G$ and their
classification is undertaken in section \ref{Sec_Twisted_Conj_Classes}.
As in the classical theory, there are twisted versions of the Weyl
group, of the conjugation map $G/T\times T\rightarrow G$ and of a
fundamental alcove parametrizing the orbits of $\mathrm{Ad}^{\kappa}$.
The main results pertaining to these objects in \cite{[Mo00],[Wen01],[MW04]}
are re-derived after establishing certain key properties of the $\kappa$-twisted
normalizer $N_{G}^{\kappa}(T^{\kappa})=\{g\in G\mbox{ }|\mbox{ }\mathrm{Ad}_{g}^{\kappa}(T^{\kappa})=T^{\kappa}\}$.

The next section deals with Fuchs, Schellekens and Schweigert's twining
characters $\tilde{\chi}^{(\kappa)}:G\rightarrow\mathbb{C}$, which
are the natural twisted generalization of the usual characters. Indeed,
these class functions satisfy orthogonality relations, give a Fourier
basis of the space of $\mathrm{Ad}^{\kappa}(G)$-invariant $L^{2}$-functions,
and generate $\widetilde{R}^{(\kappa)}(G)$ and $\widetilde{R}_{k}^{(\kappa)}(G)$,
the twisted counterparts to the representation ring and the fusion
ring at level $k\in\mathbb{Z}_{>0}$. Most notably, if $G_{(\kappa)}$
is the orbit Lie group (simply connected with root system $\mathfrak{R}_{(\kappa)}$),
we have that $\widetilde{R}_{k}^{(\kappa)}(G)\simeq R_{k}(G_{(\kappa)})$.
This is a recent result of Hong \cite{[Hon17a],[Hon17b]} that follows
from Jantzen's character formula \cite{[Jant73], [FSS96],[Wen01]},
and that is recovered here with a different set of techniques.

\subsection*{Acknowledgements}

The author is indebted to his thesis advisor, Prof. Meinrenken, for
introducing him to this subject, for his guidance, and for his comments
on this work. We also thank Prof. Jeffrey and Dr. Loizides for advice,
discussion and comments. This research was partially supported by
FRQNT and OGS doctoral scholarships.

\section{\textbf{Diagram Automorphisms and Associated Root Systems \label{Sec_Diag_Aut_Root_Systs}}}

The aim of this section is to set-up the notations used in this work,
and to present some reminders on the root systems that one can associate
to an indecomposable root system $\mathfrak{R}$ and a Dynkin diagram
automorphism $\kappa\in\mathrm{Aut}(\Pi)$, where $\Pi\subseteq\mathfrak{R}$
is a choice of simple roots. For the sake of brevity, we do not give
detailed proofs here, and refer the reader to \cite[Ch.1]{[Mo00]}.

\subsection{Notation\label{Subsection_Basics_Notation}}

\subsubsection{Groups and Lie algebras}

Let $G$ be a compact, connected, simply connected and simple Lie
group with Lie algebra $\mathfrak{g}$. We consider a fixed maximal
torus $T\subseteq G$ throughout, with Lie algebra $\mathfrak{t}$
of dimension $r=\mathrm{rk}(\mathfrak{g})$, and denote by $W=N_{G}(T)/T$
the Weyl group with respect to this torus. We assume throughout that
we have an $\mathrm{Ad}(G)$-invariant non-degenerate symmetric bilinear
form $B$ on $\mathfrak{g}$, and we use the shorthand notation $\xi\cdot\zeta=B(\xi,\zeta)$
for $\xi,\zeta\in\mathfrak{g}$. 

\subsubsection{Roots and coroots}

Below, $\mathfrak{R}\subseteq\mathfrak{t}^{\ast}$ is the root system
$\mathfrak{R}(G,T)$, $\mathfrak{R}_{+}$ is a fixed choice of positive
roots with cardinality $n_{+}=|\mathfrak{R}_{+}|$, $\Pi=\{\alpha_{i}\}_{i=1}^{r}$
designates the simple roots, and $\rho=\frac{1}{2}\sum_{\alpha\in\mathfrak{R}_{+}}\alpha$
is the half-sum of positive roots of $G$. We use real roots in this
work: for $\xi\in\mathfrak{t}$ and $t=e^{\xi}\in T$, we write $t^{\alpha}=e^{2\pi i\langle\alpha,\xi\rangle}\in\mathbb{C}$
for the corresponding character. The (closed) fundamental Weyl chamber
with respect to $\mathfrak{R}_{+}$ is $\mathfrak{t}_{+}\subseteq\mathfrak{t}$,
and the closed fundamental alcove is: 
\[
\mathfrak{A}=\left\{ \xi\in\mathfrak{t}_{+}\mbox{ }\big|\mbox{ }\langle\xi,\theta\rangle\le1\right\} ,
\]
where $\theta\in\mathfrak{R}_{+}$ is the highest root of $G$. Letting
$(\cdot,\cdot)$ denote the inner product on $\mathfrak{g}^{\ast}$
induced by $B\in S^{2}\mathfrak{g}^{\ast}$ , we assume that $B$
is normalized so that for any long root $\alpha\in\mathfrak{R}$,
we have $|\!|\alpha|\!|^{2}=(\alpha,\alpha)=2$. For $\alpha\in\mathfrak{R}$,
the corresponding coroot is $\alpha^{\vee}=\frac{2}{|\!|\alpha|\!|^{2}}(\alpha,\cdot)\in\mathfrak{t}$,
and the coroots and simple coroots of $G$ are denoted by $\mathfrak{R}^{\vee}$
and $\Pi^{\vee}$.

\subsubsection{Lattices}

We denote by $\Lambda=\ker(\exp)\cap\mathfrak{t}$ the integral lattice
of $G$, and by $\Lambda^{\ast}=\mathrm{Hom}_{\mathbb{Z}}(\Lambda,\mathbb{Z})$
the character lattice of $T$. The root and coroot lattices are respectively
denoted by $Q=\mathbb{Z}[\mathfrak{R}]\subset\mathfrak{t}^{\ast}$
and $Q^{\vee}=\mathbb{Z}[\mathfrak{R}^{\vee}]\subseteq\mathfrak{t}$,
while the weight and coweight lattices are denoted by $P=\mathrm{Hom}_{\mathbb{Z}}(Q^{\vee},\mathbb{Z})\subseteq\mathfrak{t}^{\ast}$
and $P^{\vee}=\mathrm{Hom}_{\mathbb{Z}}(Q,\mathbb{Z})\subseteq\mathfrak{t}.$
For the other root systems introduced below, we will use the same
letters and adequate sub/superscripts to denote the same objects.
Since $G$ is simply connected, we sometimes use the identifications
$\Lambda=Q^{\vee}$ and $\Lambda^{\ast}=P$ in the sequel, so that
the maximal torus is $T=\mathfrak{t}/Q^{\vee}$, and the alcove is
given by $\mathfrak{A}=\mathfrak{t}/W_{\mathrm{aff}}$, where $W_{\mathrm{aff}}=\Lambda\rtimes W$
is the affine Weyl group with $\Lambda$ acting on $\mathfrak{t}$
by translations.

\subsection{Dynkin Diagram Automorphisms\label{Subsection_Basics_Dynkin-Aut}}

Considering the automorphism group $\mathrm{Aut}(G)$ of $G$, and
its normal subgroup $\mathrm{Inn}(G)\simeq G/Z(G)$ of inner automorphisms,
the group of outer automorphisms of $G$ is defined as the quotient
$\mathrm{Out}(G)=\mathrm{Aut}(G)/\mathrm{Inn}(G)$. The latter can
be identified with the automorphism group $\mathrm{Aut}(\Pi)$ of
the Dynkin diagram of $\mathfrak{g}$, as we now briefly indicate.
By a diagram automorphism $\kappa\in\mathrm{Aut}(\Pi)$, we mean a
bijection: 
\[
\kappa:\Pi\longrightarrow\Pi,\ \ \alpha_{i}\longmapsto\alpha_{\kappa(i)},
\]
that preserves the entries of the Cartan matrix associated to $\mathfrak{R}$:
\[
\langle\alpha_{\kappa(j)},\alpha_{\kappa(i)}^{\vee}\rangle=\langle\alpha_{j},\alpha_{i}^{\vee}\rangle,\ \ i,j=1,\cdots,r.
\]
Any such automorphism induces a unique element $\kappa\in\mathrm{Aut}(\mathfrak{g}_{\mathbb{C}})$
acting on the Chevalley generators $e_{\pm\alpha}\in\mathfrak{g}_{\mathbb{C}}$
associated to simple roots as: 
\begin{equation}
\kappa(e_{\pm\alpha})=e_{\pm\kappa(\alpha)},\ \ \alpha\in\Pi,\label{Eq_kappa_Presc_Chevalley}
\end{equation}
It then follows that $\kappa$ preserves the real forms $\mathfrak{g}\subseteq\mathfrak{g}_{\mathbb{C}}$
and $\mathfrak{t}\subseteq\mathfrak{t}_{\mathbb{C}}$, and exponentiates
to an automorphism of $G$ preserving $T$. Furthermore, the element
$\kappa\in\mathrm{Aut}(\mathfrak{g})$ preserves the fundamental chamber
$\mathfrak{t}_{+}$, since it permutes the positive roots $\mathfrak{R}_{+}$.
Using the identification $\mathrm{Aut}(G)\equiv\mathrm{Aut}(\mathfrak{g})$,
equation (\ref{Eq_kappa_Presc_Chevalley}) hence defines a group homomorphism
$\mathrm{Aut}(\Pi)\to\mathrm{Aut}(G)$ descending to an isomorphism
with $\mathrm{Out}(G)$. For a given $\kappa\in\mathrm{Aut}(\Pi)$,
we will also denote the corresponding elements in $\mathrm{Out}(\mathfrak{g})$
and $\mathrm{Out}(G)$ by $\kappa$.

To study the $\kappa$-twisted conjugation $\mathrm{Ad}_{g}^{\kappa}:G\rightarrow G$,
$x\mapsto gx\kappa(g^{-1})$, it is sufficient to take $\kappa\in\mathrm{Aut}(G)$
corresponding to a Dynkin diagram automorphism. To see this, notice
that for $\kappa\in\mathrm{Aut}(G)$ and $\mathrm{Ad}_{a}\in\mathrm{Inn}(G)$,
the right action map $R_{a}:G\rightarrow G$, $x\mapsto xa^{-1}$
gives an equivariant diffeomorphism from $(\mathrm{Ad}_{a}\circ\kappa)$-twisted
conjugacy classes to $\kappa$-classes, since for any $g\in G$:
\[
R_{a^{-1}}\circ\mathrm{Ad}_{g}^{(\mathrm{Ad}_{a}\circ\kappa)}=\mathrm{Ad}_{g}^{\kappa}\circ R_{a^{-1}}.
\]
We assume from now onwards that $\kappa\in\mathrm{Out}(G)$ is a non-trivial
automorphism, and therefore restrict $\mathfrak{R}$ to the cases
compiled in the table below.

\begin{table}[H]
\begin{center}%
\begin{tabular}{|c|c|c|c|c|}
\hline 
$\mathfrak{\ensuremath{R}}$ & $A_{n}\mbox{, }n\ge2$ & $D_{n+1}\mbox{, }n\ge4$ & $D_{4}$ & $E_{6}$\tabularnewline
\hline 
\hline 
$\mathrm{Out}(G)$ & $\mathbb{Z}_{2}$ & $\mathbb{Z}_{2}$ & $S_{3}$ & $\mathbb{Z}_{2}$\tabularnewline
\hline 
\end{tabular}\end{center}\caption{Non trivial $\mathrm{Out}G$}

\end{table}

Let $\kappa\in\mathrm{Out}(G)$ be a fixed diagram automorphism of
$G$, and assume that the inner-product $B$ on $\mathfrak{g}$ is
$\kappa$-invariant. The $\kappa$-fixed subgroups will be denoted
by $G^{\kappa}\subseteq G$ and $T^{\kappa}\subseteq T$, the corresponding
Lie algebras by $\mathfrak{g^{\kappa}}$ and $\mathfrak{t}^{\kappa}$
respectively, and the subgroup of elements of the Weyl group commuting
with $\kappa$ by $W^{\kappa}\subset W$. The orthogonal projection
$p:\mathfrak{g}\rightarrow\mathfrak{g}^{\kappa}$ is given by
\[
p(\xi)=\tfrac{1}{|\kappa|}\sum_{j=1}^{|\kappa|}\kappa^{j}(\xi),\ \ \xi\in\mathfrak{g},
\]
and we use the same letter to denote its restriction $\mathfrak{t}\rightarrow\mathfrak{t}^{\kappa}$
as well as the projection $\mathfrak{t}^{\ast}\rightarrow(\mathfrak{t}^{\kappa})^{\ast}$.

Lastly, we will need certain subgroups of $T$ obtained from $\kappa\in\mathrm{Out}(G)$.
Considering the $\mathrm{Ad}$- and $\kappa$-invariant inner product
$B$ on $\mathfrak{g}$, we define $\mathfrak{t}_{\kappa}:=(\mathfrak{t}^{\kappa})^{\perp}$
with respect to the restriction $B_{|\mathfrak{t}}$, as well as the
subtorus $T_{\kappa}:=\exp(\mathfrak{t}_{\kappa})\subseteq T$. The
orthogonal decomposition $\mathfrak{t}=\mathfrak{t}^{\kappa}\oplus\mathfrak{t}_{\kappa}$
yields:
\[
T=T^{\kappa}\cdot T_{\kappa}.
\]
The intersection $T^{\kappa}\cap T_{\kappa}\subseteq T$ is a finite
subgroup, and will be of particular importance in the remaining of
this paper.

\subsection{The Folded and the Orbit Root Systems\label{Subsection_Root_Systems_F-O}}

For $\alpha\in\mathfrak{R}$, let $\overset{\vee}{\alpha}=\frac{2}{|\!|\alpha|\!|^{2}}\alpha\in\mathfrak{t}^{\ast}$.
We introduce the main constructions of this section:
\begin{defn}
To the data $(\mathfrak{R},\kappa)$, we associate the following two
root systems:\begin{itemize}

\item The \textbf{folded root system}, which is the image $p(\mathfrak{R})$
under the orthogonal projection $p:\mathfrak{t}^{\ast}\longrightarrow(\mathfrak{t}^{\ast})^{\kappa}$:
\[
\mathfrak{R}_{\mathrm{F}}^{(\kappa)}=\left\{ p(\alpha)\mbox{ }|\mbox{ }\alpha\in\mathfrak{R}\right\} \subseteq(\mathfrak{t}^{\kappa})^{\ast}.
\]

\item The \textbf{orbit root system}, given by $\mathfrak{R}_{\mathrm{O}}^{(\kappa)}=\overset{\vee}{(\mathfrak{R}_{\mathrm{F}}^{(\kappa)})}$
for $\mathfrak{R}=A_{2n-1}$, $D_{n+1}\mbox{ }(n\ge3)$, and $E_{6}$.
For the case of $\mathfrak{R}=A_{2n}$ with $n>1$, we define $\mathfrak{R}_{\mathrm{O}}^{(\kappa)}$
to be the dual in $(\mathfrak{t}^{\kappa})^{\ast}$ of the $B_{n}$
subsystem of $\mathfrak{R}_{\mathrm{F}}^{(\kappa)}=BC_{n}$. For $\mathfrak{R}=A_{2}$,
$\mathfrak{R}_{\mathrm{O}}^{(\kappa)}$ is the $A_{1}$ system with
unique positive root $2(\alpha_{1}+\alpha_{2})$.

\end{itemize}
\end{defn}

Under our assumptions on $\mathfrak{R}$ and $\kappa$, the definitions
above give the table:

\begin{table}[H]
\begin{center}%
\begin{tabular}{|c|c|c|c|c|c|c|}
\hline 
$\mathfrak{R}$ & $A_{2n-1}$ & $A_{2n}$, $n>1$ & $A_{2}$ & $D_{n+1}\mbox{, }|\kappa|=2$ & $D_{4}\mbox{, }|\kappa|=3$ & $E_{6}$\tabularnewline
\hline 
\hline 
$\mathfrak{R}_{\mathrm{F}}^{(\kappa)}$ & $C_{n}$ & $BC_{n}$ & $A_{1}\sqcup A_{1}$ & $B_{n}$ & $G_{2}$ & $F_{4}$\tabularnewline
\hline 
$\mathfrak{R}_{\mathrm{O}}^{(\kappa)}$ & $B_{n}$ & $C_{n}$ & $A_{1}$ & $C_{n}$ & $G_{2}$ & $F_{4}$\tabularnewline
\hline 
\end{tabular}\end{center}

\caption{Folded and orbit root systems}
\label{Table_RootSyst}
\end{table}

\begin{rem}
\noindent\begin{enumerate}

\item For the proof that $\mathfrak{R}_{\mathrm{F}}^{(\kappa)}=p(\mathfrak{R})$
is indeed a root system with Weyl group $W^{\kappa}$, see \cite[Prop.13.2.2]{[CaSGLT72]}.
By definition, $W^{\kappa}$ is also the Weyl group of $\mathfrak{R}_{\mathrm{O}}^{(\kappa)}$.

\item The case $\mathfrak{R}=A_{2n}$ is the only one where $\mathfrak{R}_{\mathrm{F}}^{(\kappa)}$
is decomposable. When $n=1$, $\mathfrak{R}_{\mathrm{F}}^{(\kappa)}$
is the disjoint union of two copies of $A_{1}$ (see next lemma).
When $n>1$, there are 3 root lengths in $\mathfrak{R}_{\mathrm{F}}^{(\kappa)}=BC_{n}$,
with the long and intermediate roots constituting a $C_{n}$ subsystem,
and the short and intermediate roots forming a $B_{n}$ subsystem.

\item The terminology ``orbit root system'' for $\mathfrak{R}_{\mathrm{O}}^{(\kappa)}$
follows Fuchs-Schellekens-Schweigert \cite{[FSS96]} and Hong \cite{[Hon17a],[Hon17b]}.
This is the root system denoted by $\mathfrak{R}'$ in \cite{[Mo00]}
and by $\mathfrak{R}^{1}$ in \cite{[Wen01]}, which is obtained by
rescaling the elements of $\mathfrak{R}_{\mathrm{F}}^{(\kappa)}$.

\end{enumerate}

In order to lighten the notation, we drop the $(\kappa)$ superscripts
for the remaining of this subsection, bearing in mind that we are
discussing objects associated to a fixed automorphism $\kappa\in\mathrm{Out}(G)$.
\end{rem}

\begin{lem}
\label{Lemma_Folded_Orbit_Simple_Roots} Let $(\mathfrak{R},\kappa)$
be as in table \ref{Table_RootSyst}, and denote by $\Pi^{\kappa}\subset\Pi$
the $\kappa$-invariant simple roots.\begin{itemize}

\item[(i)] For $\mathfrak{R}\ne A_{2n}$, the simple roots of $\mathfrak{R}_{\mathrm{F}}$
are given by $\Pi_{\mathrm{F}}=\{p(\alpha)\mbox{ }|\mbox{ }\alpha\in\Pi\}$,
and those of $\mathfrak{R}_{\mathrm{O}}$ by:
\[
\Pi_{\mathrm{O}}=\left\{ \alpha\mbox{ }|\mbox{ }\alpha\in\Pi^{\kappa}\right\} \sqcup\left\{ |\kappa|p(\alpha)\mbox{ }|\mbox{ }\alpha\in\Pi\smallsetminus\Pi^{\kappa}\right\} .
\]

\item[(ii)] For $\mathfrak{R}=A_{2n}$ with $n>1$, the simple roots
of the $B_{n}$ and $C_{n}$ subsystems of $\mathfrak{R}_{\mathrm{F}}$
are:
\begin{eqnarray*}
\Pi_{\mathrm{F},B_{n}} & = & \left\{ p(\alpha_{1}),\cdots,p(\alpha_{n-1}),p(\alpha_{n})\right\} ,\\
\Pi_{\mathrm{F},C_{n}} & = & \left\{ p(\alpha_{1}),\cdots,p(\alpha_{n-1}),2p(\alpha_{n})\right\} ,
\end{eqnarray*}
and the simple roots of $\mathfrak{R}_{\mathrm{O}}$ are given by:
\[
\Pi_{\mathrm{O}}=\left\{ 2p(\alpha_{1}),\cdots,2p(\alpha_{n-1}),4p(\alpha_{n})\right\} .
\]

\item[(iii)] For $\mathfrak{R}=A_{2}$, we have $\mathfrak{R}_{\mathrm{F}}=\{\pm\frac{1}{2}(\alpha_{1}+\alpha_{2}),\pm(\alpha_{1}+\alpha_{2})\}$
and $\mathfrak{R}_{\mathrm{O}}=\{\pm2(\alpha_{1}+\alpha_{2})\}$.

\end{itemize}
\end{lem}

\noindent\emph{Outline of proof.} Let $\mathfrak{R}^{\kappa}\subseteq\mathfrak{R}$
denote the subset $\kappa$-invariant roots.\begin{itemize}

\item[(i)] When $\mathfrak{R}=A_{2n-1}$, $D_{n+1}$ and $E_{6}$,
a naive identification of the nodes in a $\kappa$-orbit of the diagram
of $\mathfrak{R}$ gives the diagram of $\mathfrak{R}_{\mathrm{F}}$,
with simple roots given by $\Pi_{\mathrm{F}}=\{p(\alpha)\mbox{ }|\mbox{ }\alpha\in\Pi\}$,
and lengths of the other elements of $\mathfrak{R}_{\mathrm{F}}$
determined by whether or not $\kappa(\alpha)=\alpha$. The definition
of the orbit root system says that:
\[
\mathfrak{R}_{\mathrm{O}}=\left\{ \tfrac{2p(\alpha)}{\left(p(\alpha),p(\alpha)\right)}\mbox{ }\Big|\mbox{ }\alpha\in\mathfrak{R}\right\} \subseteq(\mathfrak{t}^{\kappa})^{\ast},
\]
which can be re-expressed as:
\[
\mathfrak{R}_{\mathrm{O}}=\left\{ \alpha\ |\ \alpha\in\mathfrak{R}^{\kappa}\right\} \cup\left\{ |\kappa|p(\alpha)\ |\ \alpha\in\mathfrak{R}\smallsetminus\mathfrak{R}^{\kappa}\right\} .
\]
For the $\kappa$-fixed $\alpha\in\mathfrak{R}$, we have $p(\alpha)=\alpha$
with $|\!|\alpha|\!|^{2}=2$, while for the elements satisfying $\kappa(\alpha)\ne\alpha$,
we have that $|\!|p(\alpha)|\!|^{2}=\frac{1}{|\kappa|}|\!|\alpha|\!|^{2}$.
The set $\Pi_{\mathrm{O}}$ in the statement gives a set of simple
roots for $\mathfrak{R}_{\mathrm{O}}$.

\item[(ii)] In the case of $\mathfrak{R}=A_{2n}$ with $n>1$, we
have a partition $\mathfrak{R}=\mathfrak{\mathfrak{R}}^{\kappa}\sqcup\mathfrak{\mathfrak{R}}^{(d)}\sqcup\mathfrak{\mathfrak{R}}^{(q)}$,
where:
\begin{eqnarray*}
\mathfrak{\mathfrak{R}}^{(d)} & = & \left\{ \alpha\in\mathfrak{R}\mbox{ }\big|\mbox{ }\kappa(\alpha)\ne\alpha\mbox{ and }(\kappa(\alpha),\alpha)=0\right\} ,\\
\mathfrak{\mathfrak{R}}^{(q)} & = & \left\{ \alpha\in\mathfrak{R}\mbox{ }\big|\mbox{ }\kappa(\alpha)\ne\alpha\mbox{ and }(\kappa(\alpha),\alpha)=-1\right\} .
\end{eqnarray*}
Projecting onto $(\mathfrak{t}^{\kappa})^{\ast}$, the elements of
$\mathfrak{R}^{\kappa}=p(\mathfrak{R}^{\kappa})$ have length $2$,
those of $p(\mathfrak{\mathfrak{R}}^{(d)})$ have length $1$, and
those of $p(\mathfrak{\mathfrak{R}}^{(q)})$ have length $\frac{1}{2}$.
In $\mathfrak{R}_{\mathrm{F}}$, we see from explicit verifications
that $p(\mathfrak{R}^{\kappa}\sqcup\mathfrak{\mathfrak{R}}^{(d)})$
is a $C_{n}$ subsystem with simple roots $\Pi_{\mathrm{F},C_{n}}$,
while $p(\mathfrak{\mathfrak{R}}^{(q)}\sqcup\mathfrak{\mathfrak{R}}^{(d)})$
constitutes a $B_{n}$ subsystem with simple roots $\Pi_{\mathrm{F},B_{n}}$.
Here, the orbit root system can be described as:
\[
\mathfrak{R}_{\mathrm{O}}=\left\{ 2\alpha\ |\ \alpha\in\mathfrak{R}^{\kappa}\right\} \cup\left\{ 2p(\alpha)\ \big|\ \alpha\in\mathfrak{R}^{(d)}\right\} ,
\]
with $2\mathfrak{R}^{\kappa}=4p(\mathfrak{R}^{(q)})$, and a set of
simple roots is given by: 
\[
\Pi_{\mathrm{O}}=\left\{ \tfrac{2}{|\!|\alpha|\!|^{2}}\alpha\mbox{ }\Big|\mbox{ }\alpha\in\Pi_{\mathrm{F},B_{n}}\right\} .
\]

\item[(iii)] Follows from the definition. \qed

\end{itemize}

\hfill

We want to realize $\mathfrak{R}_{\mathrm{F}}$ and $\mathfrak{R}_{\mathrm{O}}$
as the root systems of Lie groups related to $G$. The next proposition
determines the (co)root and (co)weight lattices of these two root
systems in terms of those of $\mathfrak{R}$ (see also \cite[Lemmas 1.1, 1.4]{[Mo00]}),
which will yield the root data needed. Before doing so, recall that
we introduced the subtorus $T_{\kappa}\subseteq T$ at the end of
\ref{Subsection_Basics_Notation}, which gives rise to the finite
subgroup $T^{\kappa}\cap T_{\kappa}\subseteq T^{\kappa}$. We define
the lattice:
\begin{equation}
\Lambda_{(\kappa)}=\exp_{\mathfrak{t}^{\kappa}}^{-1}(T^{\kappa}\cap T_{\kappa})\subseteq\mathfrak{t}^{\kappa}.\label{Eq_Orbit_Lattice}
\end{equation}
We now have the following facts:
\begin{prop}
\label{Prop_Lattices}With the notations of subsection \ref{Subsection_Basics_Notation},
we have that:\begin{enumerate}

\item If $\mathfrak{R}=A_{2n-1}$, $D_{n}$ and $E_{6}$, the lattices
of $\mathfrak{R}_{\mathrm{F}}$ satisfy the inclusions:
\[
Q_{\mathrm{F}}^{\vee}=\Lambda^{\kappa}\subseteq P_{\mathrm{F}}^{\vee}=(P^{\vee})^{\kappa},
\]
\[
Q_{\mathrm{F}}=p(Q)\subseteq p(\Lambda^{\ast})=P_{\mathrm{F}}.
\]
For the orbit root system $\mathfrak{R}_{\mathrm{O}}$, we have that
$\Lambda_{(\kappa)}=p(\Lambda)$, along with the inclusions:
\[
Q_{\mathrm{O}}^{\vee}=p(\Lambda)\subseteq P_{\mathrm{O}}^{\vee}=p(P^{\vee})
\]
\[
Q_{\mathrm{O}}=Q^{\kappa}\subseteq(\Lambda^{\ast})^{\kappa}=P_{\mathrm{O}}.
\]

\item If $\mathfrak{R}=A_{2n}$, let $P_{B_{n}}$, $Q_{B_{n}}$,
$P_{B_{n}}^{\vee}$ and $Q_{B_{n}}^{\vee}$ denote the lattices associated
to the $B_{n}$ subsystem of $\mathfrak{R}_{\mathrm{F}}=BC_{n}$.
Then:
\[
Q_{B_{n}}^{\vee}\subseteq\Lambda^{\kappa}\subseteq P_{B_{n}}^{\vee}=(P^{\vee})^{\kappa},
\]
\[
Q_{B_{n}}=p(Q)\subseteq p(\Lambda^{\ast})\subseteq P_{B_{n}},
\]
where $\left[\Lambda^{\kappa}:Q_{B_{n}}^{\vee}\right]=\left[P_{B_{n}}:p(\Lambda^{\ast})\right]=2$.
For the orbit root system, we have that $\Lambda_{(\kappa)}=p(\Lambda)$,
along with the inclusions:
\[
Q_{\mathrm{O}}^{\vee}=p(\Lambda)\subseteq p(P^{\vee})\subseteq P_{\mathrm{O}}^{\vee},
\]
\[
Q_{\mathrm{O}}\subseteq Q^{\kappa}\subseteq(\Lambda^{\ast})^{\kappa}=P_{\mathrm{O}},
\]
with $\left[P_{\mathrm{O}}^{\vee}:p(P^{\vee})\right]=\left[Q^{\kappa}:Q_{\mathrm{O}}\right]=2$.

\end{enumerate}
\end{prop}

\noindent\emph{Outline of proof.} Firstly, the simple roots in lemma
\ref{Lemma_Folded_Orbit_Simple_Roots} determine all the simple coroots
and simple (co)weights involved, which allows to check the statements
above explicitly. We illustrate for the case of $\mathfrak{R}\ne A_{2n}$.
Let $I=\{1,\cdots,r\}$ be the indexing set for the simple roots $\Pi=\{\alpha_{i}\}_{i\in I}$,
and partition it as:
\[
I=I_{0}\sqcup I_{1}\sqcup\kappa(I_{1})\sqcup\cdots\sqcup\kappa^{|\kappa|-1}(I_{1}),
\]
with $I_{0}\subseteq I$ the $\kappa$-fixed indices. The $\kappa$-orbits
in $I$ are then partioned as $\bar{I}=\bar{I}_{0}\sqcup\bar{I}_{1}$
(with $I_{j}=\bar{I}_{j}$, $j=0,1$). Next, let: 
\[
\Pi^{\vee}=\{\alpha_{i}^{\vee}\}_{i\in I}\mbox{, }\mathfrak{F}=\{\varpi_{i}\}_{i\in I}\mbox{, and }\mathfrak{F}^{\vee}=\{\varpi_{i}^{\vee}\}_{i\in I}
\]
denote the simple coroots, the fundamental weights and fundamental
coweights of $\mathfrak{R}$ respectively. Using lemma \ref{Lemma_Folded_Orbit_Simple_Roots},
one checks that for the folded system $\mathfrak{R}_{\mathrm{F}}$:
\begin{eqnarray*}
\Pi_{\mathrm{F}}^{\vee} & = & \{\alpha_{\bar{i}}^{\vee}\}_{\bar{i}\in\bar{I}_{0}}\sqcup\{|\kappa|p(\alpha_{\bar{i}}^{\vee})\}_{\bar{i}\in\bar{I}_{1}},\\
\mathfrak{F}_{\mathrm{F}} & = & \{p(\varpi_{\bar{i}})\}_{\bar{i}\in\bar{I}},\\
\mathfrak{F}_{\mathrm{F}}^{\vee} & = & \{\varpi_{\bar{i}}^{\vee}\}_{\bar{i}\in\bar{I}_{0}}\sqcup\{|\kappa|p(\varpi_{\bar{i}}^{\vee})\}_{\bar{i}\in\bar{I}_{1}},
\end{eqnarray*}
and then use these bases to show that:

\[
\begin{array}{ccccccc}
Q_{\mathrm{F}} & = & \mathbb{Z}\left[\Pi_{\mathrm{F}}\right] & = & p(Q),\\
Q_{\mathrm{F}}^{\vee} & = & \mathbb{Z}\left[\Pi_{\mathrm{F}}^{\vee}\right] & = & \Lambda\cap\mathfrak{t}^{\kappa} & = & \Lambda^{\kappa},\\
P_{\mathrm{F}} & = & \mathbb{Z}\left[\mathfrak{F}_{\mathrm{F}}\right] & = & p(\Lambda^{\ast}) & = & \mathrm{Hom}(Q_{\mathrm{F}}^{\vee},\mathbb{Z}),\\
P_{\mathrm{F}}^{\vee} & = & \mathbb{Z}\left[\mathfrak{F}_{\mathrm{F}}^{\vee}\right] & = & P^{\vee}\cap\mathfrak{t}^{\kappa} & = & \mathrm{Hom}(Q_{\mathrm{F}},\mathbb{Z}).
\end{array}
\]
For the orbit root system $\mathfrak{R}_{\mathrm{O}}$, one checks
with $\Pi_{\mathrm{O}}=\{\alpha_{\bar{i}}\}_{\bar{i}\in\bar{I}_{0}}\sqcup\{|\kappa|p(\alpha_{\bar{i}})\}_{\bar{i}\in\bar{I}_{1}}$
that:
\begin{eqnarray*}
\Pi_{\mathrm{O}}^{\vee} & = & \{p(\alpha_{\bar{i}}^{\vee})\}_{\bar{i}\in\bar{I}},\\
\mathfrak{F}_{\mathrm{O}} & = & \{\varpi_{\bar{i}}\}_{\bar{i}\in\bar{I}_{0}}\sqcup\{|\kappa|p(\varpi_{\bar{i}})\}_{\bar{i}\in\bar{I}_{1}},\\
\mathfrak{F}_{\mathrm{O}}^{\vee} & = & \{p(\varpi_{\bar{i}}^{\vee})\}_{\bar{i}\in\bar{I}},
\end{eqnarray*}
and then that:

\[
\begin{array}{ccccccccc}
Q_{\mathrm{O}} & = & \mathbb{Z}\left[\Pi_{\mathrm{O}}\right] & = & Q\cap(\mathfrak{t}^{\kappa})^{\ast} & = & Q^{\kappa},\\
Q_{\mathrm{O}}^{\vee} & = & \mathbb{Z}\left[\Pi_{\mathrm{O}}^{\vee}\right] & = & p(\Lambda) & = & \Lambda_{(\kappa)},\\
P_{\mathrm{O}} & = & \mathbb{Z}\left[\mathfrak{F}_{\mathrm{O}}\right] & = & \Lambda^{\ast}\cap(\mathfrak{t}^{\kappa})^{\ast} & = & (\Lambda^{\ast})^{\kappa} & = & \mathrm{Hom}(Q_{\mathrm{O}}^{\vee},\mathbb{Z}),\\
P_{\mathrm{O}}^{\vee} & = & \mathbb{Z}\left[\mathfrak{F}_{\mathrm{O}}^{\vee}\right] & = & p(P^{\vee}) & = & \mathrm{Hom}(Q_{\mathrm{O}},\mathbb{Z}).
\end{array}
\]
In the case of $\mathfrak{R}=A_{2n}$, the verifications are done
using the simple roots $\Pi_{\mathrm{F},B_{n}}$ of the $B_{n}$ subsystem
of $\mathfrak{R}_{\mathrm{F}}=BC_{n}$ (dual to $\mathfrak{R}_{\mathrm{O}}$).
Also, it is easily deduced that:
\[
\left[\Lambda^{\kappa}:Q_{B_{n}}^{\vee}\right]=\left[P_{B_{n}}:p(\Lambda^{\ast})\right]=\left[P_{\mathrm{O}}^{\vee}:p(P^{\vee})\right]=\left[Q^{\kappa}:Q_{\mathrm{O}}\right]=2,
\]
from the expressions of the simple (co)roots and fundamental (co)weights.
\qed
\begin{rem}
The global idea underlying the computations above is as follows. Given
$\Lambda=Q^{\vee}\subset\mathfrak{t}$ and $\Lambda^{\ast}=P\subset\mathfrak{t}^{\ast}$,
we roughly have two non-equivalent operations on these groups that
yield lattices in $\mathfrak{t}^{\kappa}$ and $(\mathfrak{t}^{\kappa})^{\ast}$:
we can either intersect with $\mathfrak{t}^{\kappa}$ (resp. $(\mathfrak{t}^{\kappa})^{\ast}$)
to get $\kappa$-invariants, or project onto $\mathfrak{t}^{\kappa}$
(resp. $(\mathfrak{t}^{\kappa})^{\ast}$) using $p$, and one operation
is dual to the other over $\mathbb{Z}$.
\end{rem}

The maximal tori of the compact connected groups associated to $\mathfrak{R}_{\mathrm{F}}$
and $\mathfrak{R}_{\mathrm{O}}$ are given by:
\[
T_{\mathrm{F}}=\mathfrak{t}^{\kappa}/\Lambda^{\kappa}=T^{\kappa}\mbox{; }T_{\mathrm{O}}=\mathfrak{t}^{\kappa}/\Lambda_{(\kappa)}=T^{\kappa}/(T^{\kappa}\cap T_{\kappa}).
\]
On the one hand, the group $T_{\mathrm{O}}$ is the maximal torus
of a Lie group of type $\mathfrak{R}_{\mathrm{O}}$ with fundamental
group isomorphic to $\Lambda_{(\kappa)}/Q_{\mathrm{O}}^{\vee}=\{1\}$.
This motivates our next definition:
\begin{defn}
With $G$ and $\kappa\in\mathrm{Out}(G)$ as in the last paragraph,
we define the \textbf{orbit Lie group} $G_{(\kappa)}$ to be the \textit{simply
connected} Lie group with maximal torus $T_{(\kappa)}=T^{\kappa}/(T^{\kappa}\cap T_{\kappa})$
and root system $\mathfrak{R}_{(\kappa)}:=\mathfrak{R}_{\mathrm{O}}$.
\end{defn}

\noindent On the other hand, the torus $T_{\mathrm{F}}=T^{\kappa}$
is the maximal torus of the fixed-point subgroup $G^{\kappa}\subseteq G$,
which is always connected for $G$ simply connected \cite[Cor.3.15]{[DKLG00]}.
We may state:
\begin{prop}
With the notations of this subsection:\begin{enumerate}

\item If $\mathfrak{R}=A_{2n}$, then the root system of $G^{\kappa}$
coincides with the $B_{n}$ subsystem of $\mathfrak{R}_{\mathrm{F}}=BC_{n}$
and $\pi_{1}(G^{\kappa})\simeq\mathbb{Z}_{2}$. In the remaining cases,
the group $G^{\kappa}$ has root system $\mathfrak{R}_{\mathrm{F}}$
and is simply connected.

\item The orbit Lie group $G_{(\kappa)}$ is related to $G$ by the
isomorphism $^{L}G_{(\kappa)}\simeq(\phantom{}^{L}G)_{0}^{\kappa}$,
where $\phantom{}^{L}G$ denotes the Langlands dual of $G$.

\end{enumerate}
\end{prop}

\begin{rem}
The first statement above is proved in the case of $G$ not necessarily
simply connected in \cite[\S\S 13.3,14.4]{[CaSGLT72]} (see also \cite[Prop.1.1]{[Mo00]}),
and $\pi_{1}(G^{\kappa})$ can be deduced from the proposition dealing
with lattices. The second statement is discussed in \cite{[KLP09]}
for complex groups (see also the more general \cite[Cor.1]{[Spr06]},
and \cite[Rk.2.3]{[HS15]}), and follows from the definition of the
Langlands dual group along with the fact that the identity component
$(\phantom{}^{L}G)_{0}^{\kappa}$ is of the same type as $G^{\kappa}$.
\end{rem}

We close this section with a lemma that will be used for several formulas,
including the description of the fundamental alcove for the $\kappa$-twisted
adjoint action on $G$ (see \cite[\S 2]{[Hon17a]} for a more extensive
treatment). 
\begin{lem}
\label{Lem_Special_Roots} \noindent\begin{itemize}

\item[(a)] Let $\theta$ denote the highest root of $\mathfrak{R}$.
The highest root $\theta_{(\kappa),\mathrm{l}}$ and the highest short
root $\theta_{(\kappa),\mathrm{s}}$ of $\mathfrak{R}_{(\kappa)}$are
given by:
\[
\theta_{(\kappa),\mathrm{l}}=\begin{cases}
|\kappa|\theta_{\mathrm{s}}^{\kappa}, & \mathfrak{R}\ne A_{2n};\\
4\theta_{\mathrm{s}}^{\kappa}=2\theta, & \mathfrak{R}=A_{2n}\mbox{, }n>1;\\
2(\alpha_{1}+\alpha_{2}), & \mathfrak{R}=A_{2};
\end{cases}
\]
\[
\theta_{(\kappa),\mathrm{s}}=\theta=\begin{cases}
\theta_{\mathrm{l}}^{\kappa}, & \mathfrak{R}\ne A_{2n};\\
4\theta_{\mathrm{l}}^{\kappa}, & \mathfrak{R}=A_{2n}\mbox{, }n>1;
\end{cases}
\]
where $\theta_{\mathrm{l}}^{\kappa}$ and $\mathrm{\theta}_{\mathrm{s}}^{\kappa}$
denote the highest root and the highest short root of $\mathfrak{R}_{\mathrm{F}}$
(resp. of $\mathfrak{R}_{B_{n}}$) for $\mathfrak{R}\ne A_{2n}$ (resp.
$\mathfrak{R}=A_{2n}$).

\item[(b)] The half-sums of positive roots of $\mathfrak{R}$ and
$\mathfrak{R}_{(\kappa)}$ are equal:
\[
\rho=\tfrac{1}{2}\sum_{\alpha\in\mathfrak{R}_{+}}\alpha=\tfrac{1}{2}\sum_{\alpha\in\mathfrak{R}_{(\kappa)+}}\alpha.
\]
\end{itemize}
\end{lem}

\begin{proof}
Part (a) can be checked using the tables on affine root systems in
\cite{[BouLieVol2], CaLAFAT05, [KacIDLA90]}. To see that $\rho$
is also the half-sum of positive roots of $\mathfrak{R}_{(\kappa)}$,
it suffices to write it as the sum of fundamental weights $\{\varpi_{i}\}_{i=1}^{r}\subset\mathfrak{t}^{\ast}$:
\[
\rho=\sum_{i=1}^{r}\varpi_{i}=\sum_{i=\kappa(i)}\varpi_{i}+\sum_{i\ne\kappa(i)}|\kappa|p(\varpi_{i}).
\]
The forms $\{\varpi_{i}|i=\kappa(i)\}\sqcup\{|\kappa|p(\varpi_{i})|i\ne\kappa(i)\}$
are the fundamental weights of $\mathfrak{R}_{(\kappa)}$ used in
the proof of proposition \ref{Prop_Lattices}.
\end{proof}
\begin{rem}
Note that since $\rho\in(\mathfrak{t}^{\kappa})^{\ast}$, there exists
a regular element in $T^{\kappa}$ under the adjoint action $\mathrm{Ad}$.
This fact will be used several times in the next section.
\end{rem}

\section{\textbf{Outer Weyl Group and Conjugacy Classes} \label{Sec_Twisted_Conj_Classes}}

In this section, we look at the generalizations to twisted conjugation
of several results in Weyl's classical theory of compact groups. In
the first subsection, we study the substitute to the Weyl group $W^{(\kappa)}$.
The main result of the second subsection is proposition \ref{Prop_Twisted_Alcove},
and says that the $\kappa$-twisted conjugacy classes in $G$ are
parametrized by the fundamental alcove of the orbit group $G_{(\kappa)}$.

\subsection{The Structure of the Outer Weyl Group}

At the end of subsection \ref{Subsection_Basics_Dynkin-Aut}, we introduced
the group $T_{\kappa}$ with Lie algebra $\mathfrak{t}_{\kappa}=(\mathfrak{t}^{\kappa})^{\perp}\subseteq\mathfrak{t}$,
along with the finite subgroup $T^{\kappa}\cap T_{\kappa}$. By introducing
the map $\phi_{\kappa}:T\rightarrow T$, $t\mapsto t\kappa(t^{-1})$,
we can describe the subtori of $T$ associated to $\kappa$ as:
\[
T^{\kappa}=\ker(\phi_{\kappa}),\ \ T_{\kappa}=\mathrm{im}(\phi_{\kappa}),
\]
with a similar description for the Lie algebras $\mathfrak{t}^{\kappa}$
and $\mathfrak{t}_{\kappa}$. The map $\phi_{\kappa}:T\to T$ will
be convenient for the study of the outer Weyl group.
\begin{defn}
The \textbf{$\bm{\kappa}$}-\textbf{twisted normalizer} of $T^{\kappa}$
in $G$ is defined as:
\[
N_{G}^{\kappa}(T^{\kappa}):=\left\{ g\in G\mbox{ }\big|\mbox{ }\mathrm{Ad}_{g}^{\kappa}(T^{\kappa})=T^{\kappa}\right\} ,
\]
and we define the \textbf{outer Weyl group} as the quotient $W^{(\kappa)}:=N_{G}^{\kappa}(T^{\kappa})/T^{\kappa}$. 
\end{defn}

\noindent Recall that $W^{\kappa}\subseteq W$ is the subgroup of
Weyl group elements commuting with $\kappa$. Since $N_{G^{\kappa}}(T^{\kappa})\subseteq N_{G}^{\kappa}(T^{\kappa})$,
we have $W^{\kappa}\subseteq W^{(\kappa)}$. The next lemma clarifies
how $N_{G}^{\kappa}(T^{\kappa})$ fits between $N_{G}(T)$ and $N_{G^{\kappa}}(T^{\kappa})$:
\begin{lem}
\label{Lemma_SOW_Twisted_Normalizer_Explicit_Desc} The $\kappa$-twisted
normalizer $N_{G}^{\kappa}(T^{\kappa})$ coincides with the following
subgroup of $N_{G}(T)$:
\[
N_{G}^{\kappa}(T^{\kappa})=\left\{ g\in N_{G}(T)\mbox{ }\big|\mbox{ }\mathrm{Ad}_{g}^{\kappa}(e)\in T^{\kappa}\right\} .
\]
\end{lem}

\begin{proof}
First, we prove the inclusion $N_{G}^{\kappa}(T^{\kappa})\subseteq N_{G}(T^{\kappa})$.
For $x\in T^{\kappa}$ and $g\in N_{G}^{\kappa}(T^{\kappa})$, we
see from the formula:
\begin{equation}
\mathrm{Ad}_{g}^{\kappa}(x)=\mathrm{Ad}_{g}(x)\cdot\mathrm{Ad}_{g}^{\kappa}(e),\label{Eq_Ad_Ad_kappa_rel}
\end{equation}
that $\mathrm{Ad}_{g}(x)\in T^{\kappa}$. Next, let $x\in T^{\kappa}\subseteq T$
be regular under $\mathrm{Ad}$, so that its centralizer is $Z_{G}(x)=T$.
Then $\mathrm{Ad}_{g}(x)\in T^{\kappa}$ is again regular, with centralizer:
\[
\mathrm{Ad}_{g}\left(Z_{G}(x)\right)=Z_{G}\left(\mathrm{Ad}_{g}(x)\right)=T,
\]
which shows that $\mathrm{Ad}_{g}(T)=T$, and therefore that $g\in N_{G}(T)$
for any $g\in N_{G}^{\kappa}(T^{\kappa})$.

Secondly, we define:
\[
H:=\left\{ g\in N_{G}(T)\mbox{ }\big|\mbox{ }g\kappa(g^{-1})\in T^{\kappa}\right\} ,
\]
and we claim that $H=N_{G}^{\kappa}(T^{\kappa})$. For any $x\in T^{\kappa}$
and $h\in H\subseteq N_{G}(T)$, we have that $\mathrm{Ad}_{h}(x)\in T^{\kappa}$,
since $\mathrm{Ad}_{h}(x)\in T$ and:
\begin{eqnarray*}
\kappa\left(\mathrm{Ad}_{h}(x)\right) & = & \mathrm{Ad}_{\kappa(h)}(x)=\mathrm{Ad}_{\left(h\kappa(h^{-1})\right)^{-1}}\mathrm{Ad}_{h}(x)\\
 & = & \mathrm{Ad}_{h}(x),
\end{eqnarray*}
where we used $h\kappa(h^{-1})\in T^{\kappa}$ on the first line.
Consequently, $h\in N_{G}^{\kappa}(T^{\kappa})$ by equation (\ref{Eq_Ad_Ad_kappa_rel}),
and we obtain $H\subseteq N_{G}^{\kappa}(T^{\kappa})$. The other
inclusion follows from $N_{G}^{\kappa}(T)\subseteq N_{G}(T)$ and
the fact that $\mathrm{Ad}_{g}^{\kappa}(e)\in T^{\kappa}$ for any
$g\in N_{G}^{\kappa}(T^{\kappa})$, which concludes the proof.
\end{proof}
The second lemma that we prove will be used for the characterization
of $W^{(\kappa)}$:
\begin{lem}
\label{Lemma_Representatives_kappa-fixed _W} For any $w\in W^{\kappa}\subseteq W$,
one may choose a representative $g\in N_{G}(T)$ such that $g\kappa(g^{-1})\in T^{\kappa}$.
\end{lem}

\begin{proof}
If $g_{w}\in N_{G}(T)$ is such that $g_{w}T=w\in W^{\kappa}$, then:
\[
\kappa\left(w\cdot t\right)=w\cdot\kappa(t)\mbox{, }t\in T\Longleftrightarrow\mathrm{Ad}_{\kappa(g_{w})}(t)=\mathrm{Ad}_{g_{w}}(t)\mbox{, }\forall t\in T.
\]
Taking $t\in T$ regular, we see that $g_{w}\kappa(g_{w}^{-1})$ must
be an element of $T$. On the other hand, the fact that $T=T^{\kappa}\cdot T_{\kappa}$
implies that $g_{w}\kappa(g_{w}^{-1})=t_{0}s^{-1}\kappa(s)$ for some
$t_{0}\in T^{\kappa}$ and $s\in T$, so that $(sg_{w})\kappa\left((sg_{w})^{-1}\right)=t_{0}\in T^{\kappa}$.
The representative $g\in N_{G}(T)$ of $w\in W^{\kappa}$ such that
$g\kappa(g^{-1})\in T^{\kappa}$ is then given by $g=g_{w}s'\in N_{G}(T)$,
where $s'=g_{w}^{-1}sg_{w}\in T$.
\end{proof}
In the upcoming characterization of $W^{(\kappa)}$, we will see the
appearance of the finite group $(T/T^{\kappa})^{\kappa}$, which is
in fact isomorphic to $T^{\kappa}\cap T_{\kappa}$:
\begin{lem}
\label{Lemma_SOW_Intersection_Subtori_1} Let $\pi:T\to T/T^{\kappa}$
be the canonical projection. The map: 
\[
\psi:(T/T^{\kappa})^{\kappa}\longrightarrow T^{\kappa}\cap T_{\kappa},\ \ \pi(t)\longmapsto\phi_{\kappa}(t),
\]
is an isomorphism.
\end{lem}

\begin{proof}
Let $q:N_{G}^{\kappa}(T^{\kappa})\to W^{(\kappa)}$ be the canonical
projection. It will be convenient to write $\psi$ as the composition
$\psi_{2}\circ\psi_{1}$ of two isomorphisms:
\begin{equation}
\psi_{1}:(T/T^{\kappa})^{\kappa}\longrightarrow q\left(\phi_{\kappa}^{-1}(T^{\kappa}\cap T_{\kappa})\right),\label{Eq_Iso_Finite_Tori_1}
\end{equation}
\begin{equation}
\psi_{2}:q\left(\phi_{\kappa}^{-1}(T^{\kappa}\cap T_{\kappa})\right)\longrightarrow T^{\kappa}\cap T_{\kappa},\label{Eq_Iso_Finite_Tori_2}
\end{equation}
with $q\left(\phi_{\kappa}^{-1}(T^{\kappa}\cap T_{\kappa})\right)\subseteq W^{(\kappa)}$.
On the one hand, we have $\mathrm{Ad}_{t}^{\kappa}(x)=\phi_{\kappa}(t)x\in T$
for all $x\in T^{\kappa}$ and $t\in T$, which is used to check that:
\begin{equation}
N_{G}^{\kappa}(T^{\kappa})\cap T=\phi_{\kappa}^{-1}(T^{\kappa}\cap T_{\kappa}).\label{Eq_Intersection_Subtori}
\end{equation}
The homomorphism $\pi_{1}:T\to T/T^{\kappa}$, $t\mapsto tT^{\kappa}$
is $\kappa$-equivariant, and an element $tT^{\kappa}\in T/T^{\kappa}$
is found to be $\kappa$-invariant iff $t\in\phi_{\kappa}^{-1}(T^{\kappa}\cap T_{\kappa})=N_{G}^{\kappa}(T^{\kappa})\cap T$.
The isomorphism (\ref{Eq_Iso_Finite_Tori_1}) is then given explicitly
by:
\[
\psi_{1}:(T/T^{\kappa})^{\kappa}\longrightarrow q\left(\phi_{\kappa}^{-1}(T^{\kappa}\cap T_{\kappa})\right)\mbox{, }\pi(t)\longmapsto q(t).
\]
On the other hand, the representatives of the group $q\left(\phi_{\kappa}^{-1}(T^{\kappa}\cap T_{\kappa})\right)$
in $N_{G}^{\kappa}(T^{\kappa})\cap T$ act on $T^{\kappa}$ by twisted
conjugation, which coincides with multiplication by elements of $T^{\kappa}\cap T_{\kappa}$.
The isomorphism \ref{Eq_Iso_Finite_Tori_2} is given by: 
\[
\psi_{2}:q\left(\phi_{\kappa}^{-1}(T^{\kappa}\cap T_{\kappa})\right)\longrightarrow T^{\kappa}\cap T_{\kappa}\mbox{, }tT^{\kappa}\longmapsto\phi_{\kappa}(t).\qedhere
\]
\end{proof}
The last lemma that we will need is the following:
\begin{lem}
\label{Lemma_Intersection_Tori_Expl}One has that:
\[
T^{\kappa}\cap T_{\kappa}\simeq\begin{cases}
(\mathbb{Z}_{2})^{\dim\mathfrak{t}_{\kappa}}, & |\kappa|=2;\\
\mathbb{Z}_{3}, & |\kappa|=3.
\end{cases}
\]
\end{lem}

\begin{proof}
We consider the two cases separately:\begin{itemize}

\item \textit{For $|\kappa|=2$:} Let $t=s\kappa(s^{-1})\in T^{\kappa}\cap T_{\kappa}$
be an arbitrary element, with $s\in T$. For $|\kappa|=2$, the automorphism
$\kappa$ acts as inversion on $T^{\kappa}\cap T_{\kappa}$: 
\[
\kappa(t)=\kappa(s)s^{-1}=t^{-1}=t,
\]
and it is easily checked that in fact:
\[
T^{\kappa}\cap T_{\kappa}=\left\{ t\in T_{\kappa}\mbox{ }\big|\mbox{ }t=t^{-1}\right\} ,
\]
Since $T_{\kappa}\simeq(S^{1})^{\dim\mathfrak{t}_{\kappa}}$, we see
that there are $\dim\mathfrak{t}_{\kappa}=(\dim\mathfrak{t}-\dim\mathfrak{t}^{\kappa})$
generators of order 2 for $T^{\kappa}\cap T_{\kappa}$, and the claim
follows. 

\item \textit{For $|\kappa|=3$:} In the case of $\mathfrak{R}=D_{4}$,
that $T^{\kappa}\cap T_{\kappa}\simeq\mathbb{Z}_{3}$ can be seen
more easily with a concrete description. Letting
\[
T=\left\{ (z_{1},z_{2},z_{3},z_{4})\mbox{ }\big|\mbox{ }z_{j}\in S^{1}\right\} ,
\]
the automorphism $\kappa$ acts as $\kappa(z_{1},z_{2},z_{3},z_{4})=(z_{4},z_{2},z_{1},z_{3})$.
Thus:
\begin{eqnarray*}
T^{\kappa} & = & \left\{ (a,b,a,a)\mbox{ }\big|\mbox{ }a,b\in S^{1}\right\} ,\\
T_{\kappa} & = & \left\{ (a,1,b,a^{-1}b^{-1})\mbox{ }\big|\mbox{ }a,b\in S^{1}\right\} ,
\end{eqnarray*}
and there are only 3 possibilities for $t\in T^{\kappa}\cap T_{\kappa}$:
$(e^{\pm\frac{2\pi i}{3}},1,e^{\pm\frac{2\pi i}{3}},e^{\pm\frac{2\pi i}{3}})$
and the identity.\qedhere

\end{itemize}
\end{proof}
The main result of this subsection is the following (cf. \cite[\S 2.3]{[Mo00]},
\cite[Prop.2.4]{[Wen01]}):
\begin{thm}
\label{Thm_Structure_Outer_Weyl} The outer Weyl group decomposes
as the semi-direct product:
\[
W^{(\kappa)}=(T^{\kappa}\cap T_{\kappa})\rtimes W^{\kappa}.
\]
\end{thm}

\begin{proof}
Here, we use the notation of the proof of lemma \ref{Lemma_SOW_Intersection_Subtori_1}.
We show that we have the following split exact sequence:
\[
1\longrightarrow(T/T^{\kappa})^{\kappa}\longrightarrow W^{(\kappa)}\overset{\nu}{\longrightarrow}W^{\kappa}\longrightarrow1,
\]
where: 
\[
\nu:W^{(\kappa)}\longrightarrow W^{\kappa}\mbox{, }gT^{\kappa}\longmapsto gT.
\]
We prove each claim separately. \begin{itemize}

\item\textit{ Exactness at $(T/T^{\kappa})^{\kappa}$:} If $tT^{\kappa}\in(T/T^{\kappa})^{\kappa}$
then $t\in\phi_{\kappa}^{-1}(T^{\kappa}\cap T_{\kappa})\subset N_{G}^{\kappa}(T^{\kappa})$
by equation (\ref{Eq_Intersection_Subtori}), and we have a well defined
map $\iota:(T/T^{\kappa})^{\kappa}\rightarrow W^{(\kappa)}\mbox{, }tT^{\kappa}\mapsto tT^{\kappa}$
that is clearly injective.

\item \textit{Exactness at $W^{\kappa}$:} That $\nu:W^{(\kappa)}\longrightarrow W^{\kappa}$
is a well-defined homomorphism is easily seen from its definition
above, and its surjectivity is a direct consequence of lemma \ref{Lemma_Representatives_kappa-fixed _W}.

\item \textit{Exactness at $W^{(\kappa)}$:} We show that $\ker\nu\subset(T/T^{\kappa})^{\kappa}$:
If $g\in N_{G}^{\kappa}(T^{\kappa})$ is such that $\nu(gT^{\kappa})=eT$,
then $g\in N_{G}^{\kappa}(T^{\kappa})\cap T=\phi_{\kappa}^{-1}(T^{\kappa}\cap T_{\kappa})$
and $q(g)=gT^{\kappa}\in(T/T^{\kappa})^{\kappa}$ by (\ref{Eq_Intersection_Subtori}).
Since $(T/T^{\kappa})^{\kappa}\subset\ker\nu$ is obvious, we have
$\ker\nu=(T/T^{\kappa})^{\kappa}$.

\item\textit{ Splitting:} For any $x\in T^{\kappa}$ and $gT\in W^{\kappa}$,
we have: 
\[
\mathrm{Ad}_{g}(x)=\left(g\kappa(g^{-1})\right)^{-1}\mathrm{Ad}_{g}^{\kappa}(x)\in T\mbox{, }
\]
and proceeding as in the proof of lemma \ref{Lemma_Representatives_kappa-fixed _W},
we find that for any $x\in T^{\kappa}$, there exists a $t\in T$
such that $\mathrm{Ad}_{gt}^{\kappa}(x)\in T^{\kappa}$. We thus have
an inclusion $j:W^{\kappa}\rightarrow W^{(\kappa)}\mbox{, }gT\mapsto gT^{\kappa}$
such that $\nu\circ j=\mathrm{Id}_{W^{\kappa}}$.\end{itemize}

From the above, we obtain $W^{(\kappa)}\simeq(T/T^{\kappa})^{\kappa}\rtimes W^{\kappa}$,
and the theorem follows from the isomorphism in lemma \ref{Lemma_SOW_Intersection_Subtori_1}.
\end{proof}

\subsection{Twisted Conjugacy Classes and Fundamental Alcove}

The orbit space $G/\mathrm{Ad}^{\kappa}(G)$ is controlled by the
groups $G_{(\kappa)}$, $W^{(\kappa)}$, $T^{\kappa}$ and $T_{\kappa}$.
The following lemma corresponds to \cite[Lem.2.1]{[Wen01]} (c.f.
\cite[Lem.IV.4.4]{[BtD95]}), and we present it in terms of the orbit
root system $\mathfrak{R}_{(\kappa)}$ instead of the folded root
system.
\begin{lem}
\label{Lemma_Det_Diff_Conj_Map}Define the conjugation map:
\[
\mathsf{c}:T^{\kappa}\times(G/T^{\kappa})\longrightarrow G\mbox{, }(x,gT^{\kappa})\longmapsto\mathrm{Ad}_{g}^{\kappa}(x).
\]
The determinant of the differential $d\mathsf{c}$ at the point $(t,gT^{\kappa})\in T^{\kappa}\times(G/T^{\kappa})$
is given by:
\[
\det\left(d\mathsf{c}_{|(t,gT^{\kappa})}\right)=|T^{\kappa}\cap T_{\kappa}|\cdot\big|\widetilde{\Delta}(t)\big|^{2},
\]
where $\widetilde{\Delta}:T^{\kappa}\rightarrow\mathbb{C}$ is the
function:
\[
\widetilde{\Delta}(t)=\prod_{\tilde{\alpha}\in\mathfrak{R}_{(\kappa)+}}(1-t^{-\tilde{\alpha}}).
\]
 
\end{lem}

\begin{proof}
Recall that we have a $\kappa$- and $\mathrm{Ad}$-invariant inner
product $B$ on $\mathfrak{g}$ ($\S$\ref{Subsection_Basics_Notation}).
Fixing $(t,gT^{\kappa})\in T^{\kappa}\times(G/T^{\kappa})$, the connectedness
of $G$ implies that $\mathrm{Ad}_{\kappa(g)}\in\mathrm{SO}(\mathfrak{g})$,
and in terms of the orthogonal decomposition 
\[
\left(\mathfrak{g}/\mathfrak{t}^{\kappa}\right)_{\mathbb{C}}=(\mathfrak{t}_{\kappa})_{\mathbb{C}}\oplus\bigoplus_{\alpha\in\mathfrak{R}}\mathfrak{g}_{\alpha},
\]
where $\mathfrak{g}_{\alpha}$ denotes the root space for $\alpha$,
a direct computation yields:
\[
\det\left(d\mathsf{c}_{|(t,gT^{\kappa})}\right)=\mathrm{det}_{\mathfrak{g}/\mathfrak{t}^{\kappa}}\left(\mathrm{Ad}_{t^{-1}}-\kappa\right)=\mathrm{det}_{(\mathfrak{t}_{\kappa})_{\mathbb{C}}}\left(\mathrm{Ad}_{t^{-1}}-\kappa\right)\prod_{\alpha\in\mathfrak{R}}\mathrm{det}_{\mathfrak{g}_{\alpha}}\left(\mathrm{Ad}_{t^{-1}}-\kappa\right).
\]
The determinant on $(\mathfrak{t}_{\kappa})_{\mathbb{C}}$ depends
on $|\kappa|$, and we see from an eigenspace decomposition that:
\[
\mathrm{det}_{(\mathfrak{t}_{\kappa})_{\mathbb{C}}}\left(\mathrm{Ad}_{t^{-1}}-\kappa\right)=|T^{\kappa}\cap T_{\kappa}|=\begin{cases}
2^{\dim\mathfrak{t}_{\kappa}}, & \mbox{if }|\kappa|=2;\\
3, & \mbox{if }|\kappa|=3.
\end{cases}
\]
For the remaining factors, we first note that if $e_{\alpha}\in\mathfrak{g}_{\alpha}$
is a root vector for $\alpha\in\mathfrak{R}$, then $\kappa$ acts
on it as follows (see \cite[Eq. (12.48)]{[FSLAR97]} for an explanation):
\[
\kappa(e_{\alpha})=\begin{cases}
e_{\kappa(\alpha)}, & \mbox{if }\mathfrak{R}=A_{2n-1},D_{n},E_{6};\\
(-1)^{\mathrm{ht}\alpha+1}e_{\kappa(\alpha)}, & \mbox{if }\mathfrak{R}=A_{2n}.
\end{cases}
\]
Secondly, we gather the terms in the product over $\mathfrak{R}$
by $\kappa$-orbits. Suppose for simplicity that $\mathfrak{R}\ne A_{2n}$
with $|\kappa|=2$, and that $\alpha\in\mathfrak{R}_{+}$. We obtain:
\[
\det\left(\left(\mathrm{Ad}_{t^{-1}}-\kappa\right)_{|\mathfrak{g}_{\alpha}\oplus\mathfrak{g}_{-\alpha}}\right)=|1-t^{-\alpha}|^{2}\mbox{, if }\kappa(\alpha)=\alpha,
\]
\[
\det\left(\left(\mathrm{Ad}_{t^{-1}}-\kappa\right)_{|\mathfrak{g}_{\alpha}\oplus\mathfrak{g}_{\kappa(\alpha)}\oplus\mathfrak{g}_{-\alpha}\oplus\mathfrak{g}_{-\kappa(\alpha)}}\right)=|1-t^{-|\kappa|p(\alpha)}|^{2}\mbox{, if }\kappa(\alpha)\ne\alpha.
\]
Using the description of $\mathfrak{R}_{(\kappa)}$ in the proof of
lemma \ref{Lemma_Folded_Orbit_Simple_Roots}, we have:
\[
\det\left(d\mathsf{c}_{|(t,gT^{\kappa})}\right)=|T^{\kappa}\cap T_{\kappa}|\times\prod_{\tilde{\alpha}\in\mathfrak{R}_{(\kappa)+}}|1-t^{-\tilde{\alpha}}|^{2}=|T^{\kappa}\cap T_{\kappa}|\cdot|\widetilde{\Delta}(t)|^{2}.
\]
With the appropriate modifications, the same formula is obtained for
the cases $\mathfrak{R}=A_{2n}$ and $\mathfrak{R}=D_{4}$ with $|\kappa|=3$.
\end{proof}
\noindent The second lemma needed follows from elementary considerations,
and we omit its proof.
\begin{lem}
\label{Lemma_Stabilizer_Twisted_Conj}For $x\in G$, denote by:
\[
Z_{G}^{\kappa}(x):=\left\{ g\in G\mbox{ }\big|\mbox{ }\mathrm{Ad}_{g}^{\kappa}(x)=x\right\} 
\]
the stabilizer under $\kappa$-twisted conjugation, and let $\mathfrak{z}_{x}^{\kappa}=\mathrm{Lie}\left(Z_{G}^{\kappa}(x)\right)$
be its Lie algebra. Then:\begin{enumerate}

\item For $x\in T^{\kappa}$, the group $Z_{G}^{\kappa}(x)$ is preserved
by $\kappa$ and contains $T^{\kappa}$ as a maximal torus.

\item Let $\mathcal{C}=\mathrm{Ad}^{\kappa}(G)\cdot x$. With respect
to the $\kappa$- and $\mathrm{Ad}$-invariant metric on $G$, one
has an orthogonal decomposition $T_{x}G=\mathfrak{z}_{x}^{\kappa}\oplus T_{x}\mathcal{C}$,
with $\mathfrak{z}_{x}^{\kappa}=\ker(\mathrm{Ad}_{x^{-1}}-\kappa)$,
and such that in the left trivialization of $TG$:
\[
T_{x}\mathcal{C}=\left\{ \left((\mathrm{Ad}_{x^{-1}}-\kappa)\cdot\xi\right)_{|x}^{L}\mbox{ }\big|\mbox{ }\xi\in\mathfrak{g}\right\} .
\]
(For $\xi\in\mathfrak{g}$, we use $\xi_{|g}^{L}=\frac{d}{dt}\left(g\cdot\exp(t\xi)\right)_{|t=0}$
for the left-invariant vector fields.)

\end{enumerate}
\end{lem}

The next proposition generalizes a classical result of the theory
of compact connected Lie groups, and is in essence a rephrasing of
\cite[Prop.IV.4.3]{[BtD95]} and \cite[Lem.2.1]{[Wen01]}. The analogous
result for algebraic groups is proved in \cite[\S 2.5]{[Mo00]} and
\cite[Lem.2]{[Spr06]}.
\begin{prop}
\label{Prop_Twisted_Conj_Classes}The $\kappa$-twisted conjugacy
classes in $G$ satisfy the following properties:

\begin{enumerate}

\item\label{Prop_TCC_Surj_Conj_Map} For any $g\in G$, there exist
elements $x\in T^{\kappa}$ and $a\in G$ such that $g=\mathrm{Ad}_{a}^{\kappa}(x)$.
That is, any element of $G$ is $\mathrm{Ad}^{\kappa}$-conjugate
to an element of the torus $T^{\kappa}$.

\item\label{Prop_TCC_Conjugacy_Tk} Two elements of $T^{\kappa}$
are twisted conjugate under an element of $G$ if and only if they
are twisted conjugate under an element of $N_{G}^{\kappa}(T^{\kappa})$.

\item\label{Prop_TCC_Intersection_TCC_Tk} Any $\kappa$-twisted
conjugacy class in $G$ intersects $T^{\kappa}$ in an orbit of $W^{(\kappa)}=N_{G}^{\kappa}(T^{\kappa})/T^{\kappa}$.
In particular, one has for any $x\in T^{\kappa}$ that:
\[
\mathrm{Ad}^{\kappa}(G)\cdot x\cap T^{\kappa}=W^{(\kappa)}\cdot x.
\]

\end{enumerate}
\end{prop}

\begin{proof}
\noindent \begin{itemize}

\item[(\ref{Prop_TCC_Surj_Conj_Map})] This assertion is proved by
establishing the surjectivity of the twisted conjugation map $\mathsf{c}:T^{\kappa}\times(G/T^{\kappa})\longrightarrow G$.
Let $\xi\in\mathfrak{t}^{\kappa}\cap\mathfrak{t}^{\mathrm{reg}}$
be a regular element with respect to usual conjugation, and set $x=e^{\xi}\in T^{\kappa}$.
By lemma \ref{Lemma_Det_Diff_Conj_Map}, we have $\det\left(d\mathsf{c}_{|(x,gT^{\kappa})}\right)>0$
for any $g\in G$, meaning that the element $\mathrm{Ad}_{g}^{\kappa}(x)\in G$
is a regular value of the conjugation map $\mathsf{c}$, which is
orientation preserving at $(x,gT^{\kappa})\in T^{\kappa}\times(G/T^{\kappa})$.
From the decomposition in theorem \ref{Thm_Structure_Outer_Weyl}
of $W^{(\kappa)}$, we have:
\[
\#\left(\mathsf{c}^{-1}(\mathrm{Ad}_{g}^{\kappa}(x))\right)=\#\left(W^{(\kappa)}\cdot x\right)\ge1,
\]
which shows that $\deg\mathsf{c}>0$, and that $\mathsf{c}:T^{\kappa}\times(G/T^{\kappa})\longrightarrow G$
is a surjective map.

\item[(\ref{Prop_TCC_Conjugacy_Tk})] Let $\mathcal{C}=\mathrm{Ad}^{\kappa}(G)\cdot x$
for $x\in T^{\kappa}$, and let $y=\mathrm{Ad}_{a}^{\kappa}(x)\in\mathcal{C}\cap T^{\kappa}$
for some $a\in G$. We will modify $a$ by an element $z\in Z_{G}^{\kappa}(x)$,
in such a way that $az\in N_{G}^{\kappa}(T^{\kappa})$ with $\mathrm{Ad}_{az}^{\kappa}(x)=y$,
which will prove the claim.

Consider a regular element $\xi\in\mathfrak{t}^{\kappa}$ with respect
to the usual conjugation action. Viewing $\xi$ as an element of $T_{y}(T^{\kappa})$,
the differential $(\mathrm{Ad}_{a^{-1}}^{\kappa})_{\ast|y}:T_{y}(T^{\kappa})\rightarrow T_{x}(T^{\kappa})$
sends it to $\mathrm{Ad}_{\kappa(a^{-1})}\xi\in(T_{x}\mathcal{C})^{\perp}=\mathfrak{z}_{x}^{\kappa}$
by lemma \ref{Lemma_Stabilizer_Twisted_Conj}-(2) (as $\mathfrak{t}^{\kappa}\subset\mathfrak{z}_{x}^{\kappa}$).
Now by lemma \ref{Lemma_Stabilizer_Twisted_Conj}-(1), that $T^{\kappa}$
is a maximal torus in $Z_{G}^{\kappa}(x)$ implies that there exists
an element $z\in Z_{G}^{\kappa}(x)$ such that $\mathrm{Ad}_{\kappa(z^{-1})}(\mathrm{Ad}_{\kappa(a^{-1})}\xi)\in\mathfrak{t}^{\kappa}$.
Putting $b=az\in G$, we thus have:
\[
\mathrm{Ad}_{b}^{\kappa}(x)=\mathrm{Ad}_{a}^{\kappa}(x)=y\in T^{\kappa},
\]
\[
\mathrm{Ad}_{\kappa(b^{-1})}\xi=\mathrm{Ad}_{b^{-1}}\xi\in\mathfrak{t}^{\kappa}.
\]
Since $\xi,\mathrm{Ad}_{\kappa(b^{-1})}\xi\in\mathfrak{t}$ are conjugate,
the exists some $\nu\in N_{G}(T)$ such that $\mathrm{Ad}_{(b\nu)^{-1}}\xi=\xi$,
meaning that $b\nu\in T\subset N_{G}(T)$ by regularity of $\xi$.
Thus, the first conclusion here is that:
\[
b\in N_{G}(T).
\]
On the other hand, $\mathrm{Ad}_{\kappa(b^{-1})}\xi=\mathrm{Ad}_{b^{-1}}\xi\in\mathfrak{t}^{\kappa}$
implies that $\mathrm{Ad}_{b\kappa(b^{-1})}\xi=\xi$, and by $\mathrm{Ad}$-regularity
of $\xi\in\mathfrak{t}^{\kappa}$ that $\kappa(b)=bt$ for some $t\in T$.
Proceeding then as in the proof of lemma \ref{Lemma_SOW_Twisted_Normalizer_Explicit_Desc},
we find $\kappa\left(\mathrm{Ad}_{b}(x)\right)=\mathrm{Ad}_{b}(x)\in T$
and therefore $\mathrm{Ad}_{b}(x)\in T^{\kappa}$. Since $y=\mathrm{Ad}_{b}(x)\cdot b\kappa(b^{-1})\in T^{\kappa}$,
the second conclusion is that: 
\[
b\kappa(b^{-1})\in T^{\kappa}.
\]
By lemma \ref{Lemma_SOW_Twisted_Normalizer_Explicit_Desc}, we established
that $b\in N_{G}^{\kappa}(T^{\kappa})$, which finishes the proof.

\item[(\ref{Prop_TCC_Intersection_TCC_Tk})] By the first statement,
any element $g\in G$ lies in a conjugacy class of the form $\mathcal{C}_{g}=\mathrm{Ad}^{\kappa}(G)\cdot x$
for some $x\in T^{\kappa}$. By the second statement, $y\in\mathcal{C}_{g}\cap T^{\kappa}$
if and only if there is an element $w\in W^{(\kappa)}$ such that
$y=w\cdot x$, hence the claim.\qedhere

\end{itemize}
\end{proof}
The last statement in this proposition gives an identification $G/\mathrm{Ad}^{\kappa}(G)\simeq T^{\kappa}/W^{(\kappa)}$.
To express this space of orbits as quotient of $\mathfrak{t}^{\kappa}$,
it is sufficient to determine the pre-image under $\exp_{\mathfrak{t}^{\kappa}}:\mathfrak{t}^{\kappa}\rightarrow T^{\kappa}$
of an intersection $\mathrm{Ad}^{\kappa}(G)\cdot e^{\xi}\cap T^{\kappa}=W^{(\kappa)}\cdot e^{\xi}$,
where $\xi\in\mathfrak{t}^{\kappa}$. The map $\exp_{\mathfrak{t}^{\kappa}}$
is equivariant with respect to the action of $W^{\kappa}=N_{G^{\kappa}}(T^{\kappa})/T^{\kappa}$
on $T^{\kappa}$ and $\mathfrak{t}^{\kappa}$, and using $W^{(\kappa)}\simeq(T^{\kappa}\cap T_{\kappa})\rtimes W^{\kappa}$,
we have:
\[
\exp_{\mathfrak{t}^{\kappa}}^{-1}\left(W^{(\kappa)}\cdot e^{\xi}\right)=\exp_{\mathfrak{t}^{\kappa}}^{-1}\left(W^{\kappa}\cdot\left(T^{\kappa}\cap T_{\kappa}\cdot e^{\xi}\right)\right)=W^{\kappa}\cdot\exp_{\mathfrak{t}^{\kappa}}^{-1}\left((T^{\kappa}\cap T_{\kappa})\cdot e^{\xi}\right).
\]
Recalling that we introduced $\Lambda_{(\kappa)}=\exp_{\mathfrak{t}^{\kappa}}^{-1}(T^{\kappa}\cap T_{\kappa})$
in section \ref{Subsection_Root_Systems_F-O}, this lattice has a
natural action by translations on $\mathfrak{t}^{\kappa}$, so that
$\exp_{\mathfrak{t}^{\kappa}}^{-1}\left((T^{\kappa}\cap T_{\kappa})\cdot e^{\xi}\right)=\Lambda_{(\kappa)}\cdot\xi$.
Defining the \textbf{twisted affine Weyl group} $W_{\mathrm{aff}}^{(\kappa)}:=\Lambda_{(\kappa)}\rtimes W^{\kappa}$,
the space of $\kappa$-twisted conjugacy classes is then given by:
\[
G/\mathrm{Ad}^{\kappa}(G)\simeq T^{\kappa}/W^{(\kappa)}\simeq\mathfrak{t}^{\kappa}/W_{\mathrm{aff}}^{(\kappa)}.
\]
Using the notation of \ref{Lem_Special_Roots}, we have an explicit
description of $\mathfrak{t}^{\kappa}/W_{\mathrm{aff}}^{(\kappa)}$.
\begin{prop}
\label{Prop_Twisted_Alcove}There is a unique fundamental domain containing
the origin in $\mathfrak{t}^{\kappa}$ for the action of $W_{\mathrm{aff}}^{(\kappa)}$,
which is the alcove:
\[
\mathfrak{A}^{(\kappa)}=\left\{ \xi\in\mathfrak{t}^{\kappa}\mbox{ }\big|\mbox{ }0\le\langle\tilde{\alpha},\xi\rangle\mbox{, }\forall\tilde{\alpha}\in\Pi_{(\kappa)}\mbox{; }\langle\theta_{(\kappa),\mathrm{l}},\xi\rangle\le1\right\} ,
\]
where $\theta_{(\kappa),\mathrm{l}}\in\mathfrak{R}_{(\kappa)+}$ is
the highest root of $\mathfrak{R}_{(\kappa)}$, and $\Pi_{(\kappa)}\subset(\mathfrak{t}^{\kappa})^{\ast}$
are its simple roots.
\end{prop}

\begin{example}
We illustrate this alcove in the case $G=SU(3)$ (see Fig. \ref{Fig_Twist_Alcove}).
Here $\mathfrak{R}=A_{2}$, and $\kappa\ne1$ is the automorphism
permuting the simple roots $\alpha_{1}$ and $\alpha_{2}$, the dual
subtorus $(\mathfrak{t}^{\kappa})^{\ast}\subset\mathfrak{t}^{\ast}$
is the line $\{\lambda(\alpha_{1}+\alpha_{2})\}_{\lambda\in\mathbb{R}}$,
and $\mathfrak{R}_{(\kappa)}$ has only one positive root $(\theta_{\kappa}^{2})_{\mathrm{l}}=2(\alpha_{1}+\alpha_{2})$.
In the standard coordinates on $\mathbb{R}^{2}\equiv\mathfrak{t}$,
we have $\alpha_{1}^{\vee}=\frac{\sqrt{2}}{2}(\sqrt{3},-1)$ and $\alpha_{2}^{\vee}=\sqrt{2}(0,1)$,
so that the fundamental alcove for twisted conjugation is:
\begin{eqnarray*}
\mathfrak{A}^{(\kappa)} & = & \left\{ \xi\in\mathfrak{t}^{\kappa}\mbox{ }\big|\mbox{ }0\le2\langle\alpha_{1}+\alpha_{2},\xi\rangle\le1\right\} \\
 & = & \left\{ \tfrac{1}{2}(\sqrt{3},1)t\in\mathbb{R}^{2}\mbox{ }\Big|\mbox{ }0\le t\le\tfrac{\sqrt{2}}{4}\right\} .
\end{eqnarray*}

\begin{figure}
\begin{centering}
\includegraphics[scale=0.5]{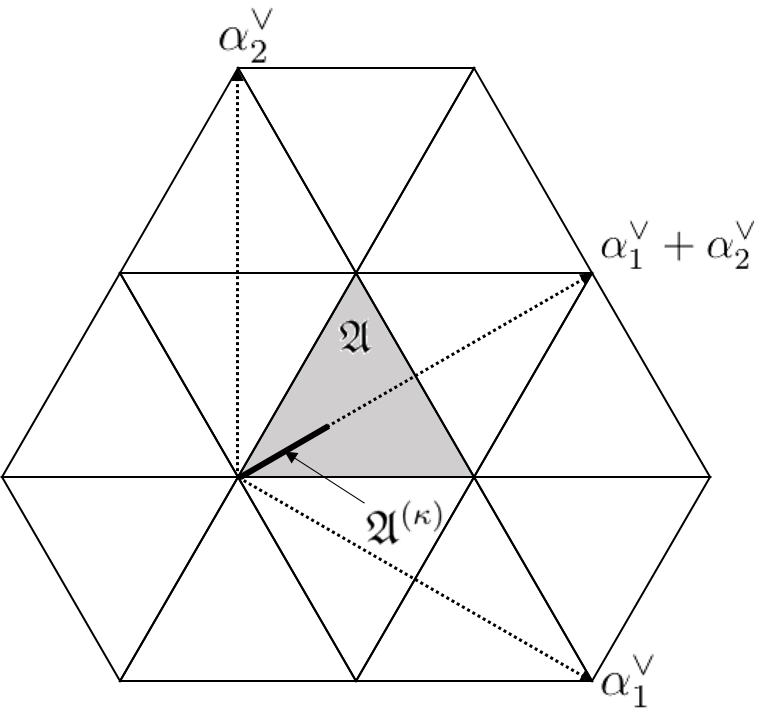}
\par\end{centering}
\caption{The alcove $\mathfrak{A}^{(\kappa)}$ for $A_{2}$}
\label{Fig_Twist_Alcove}
\end{figure}
\end{example}

\begin{rem}
The group of affine reflections $W_{\mathrm{aff}}^{(\kappa)}$ coincides
with the affine Weyl group of the untwisted affine algebra corresponding
to $\mathfrak{R}_{(\kappa)}$.
\end{rem}

For the sake of completeness, we conclude this section with Mohrdieck
and Wendt's description of the stabilizers $Z_{a}^{\kappa}=Z_{G}^{\kappa}(a)$
for $a\in T^{\kappa}$ and their Lie algebras $\mathfrak{z}_{a}^{\kappa}\subset\mathfrak{g}$\cite[Prop.4.1]{[MW04]}:
\begin{prop}
Let $\xi\in\mathfrak{A}^{(\kappa)}$, and let $\widetilde{\Pi}_{(\kappa)}=\Pi_{(\kappa)}\cup\{-\theta_{(\kappa),\mathrm{l}}\}$
label the vertices of the extended Dynkin diagram associated to $\mathfrak{R}_{(\kappa)}$.\begin{enumerate}

\item The Dynkin diagram of $(\mathfrak{z}_{e^{\xi}}^{\kappa},\mathfrak{t}^{\kappa})$
is obtained from that of $\widetilde{\Pi}_{(\kappa)}$ by deleting
the vertex $\alpha\in\Pi_{(\kappa)}$ if $\langle\alpha,\xi\rangle\ne0$,
and deleting the vertex $-\theta_{(\kappa),\mathrm{l}}$ if $\langle\theta_{(\kappa),\mathrm{l}},\xi\rangle<1$.

\item If $\mathfrak{a}\subseteq\mathfrak{g}$ denotes the subalgebra
corresponding to the diagram in (1), then: $\mathfrak{z}_{e^{\xi}}^{\kappa}=\mathfrak{t}^{\kappa}+\mathfrak{a}$.

\item The stabilizer $Z_{e^{\xi}}^{\kappa}\subset G$ is the Lie
group with the root system described in (1) and fundamental group
$\pi_{1}(Z_{e^{\xi}}^{\kappa})=\Lambda^{\kappa}/Q_{\xi}^{\vee}$,
where $Q_{\xi}^{\vee}$ is the corresponding coroot lattice.

\end{enumerate}
\end{prop}

\begin{rem}
The stabilizers $Z_{G}^{\kappa}\left(\exp_{\mathfrak{g}}(\xi)\right)\subseteq G$
and $Z_{G_{(\kappa)}}\left(\exp_{\mathfrak{g}_{(\kappa)}}(\xi)\right)\subseteq G_{(\kappa)}$
are not isomorphic for $\xi\in\mathfrak{A}^{(\kappa)}$. The simplest
example is that of $\xi=0$, for which:
\[
Z_{G}^{\kappa}(e)=G^{\kappa}\mbox{ and }Z_{G_{(\kappa)}}(e)=G_{(\kappa)}.
\]
This also shows that the $\kappa$-twisted conjugacy class corresponding
to $\xi\in\mathfrak{A}^{(\kappa)}$ is typically of higher dimension
than the orbit $\mathrm{Ad}(G_{(\kappa)})\cdot\exp_{\mathfrak{g}_{(\kappa)}}(\xi)\subseteq G_{(\kappa)}$.
\end{rem}

\section{\textbf{Twining Characters} \label{Sec_Twining_Char}}

For a compact connected Lie group $G$, let $\Lambda_{+}^{\ast}=\Lambda^{\ast}\cap\mathfrak{t}_{+}^{\ast}$
be the integral dominant weights, and denote by $\mathcal{C}^{\infty}(G)^{G}$
the ring of smooth $\mathrm{Ad}$-invariant functions on $G$. The
representation ring $R(G)$ can be realized as the subring of $\mathcal{C}^{\infty}(G)^{G}$
generated by the irreducible characters $\{\chi_{\lambda}\}_{\lambda\in\Lambda_{+}^{\ast}}$,
and as a $\mathbb{Z}$-module, one has:
\[
R(G)\simeq\mathbb{Z}[\Lambda_{+}^{\ast}].
\]
When $G$ is also simply connected, we can also define the fusion
ring at level $k\in\mathbb{N}\smallsetminus\{0\}$ solely in terms
of $\{\chi_{\lambda}\}_{\lambda\in\Lambda_{+}^{\ast}}$, by taking:
\[
R_{k}(G):=R(G)/I_{k}(G)\simeq\mathbb{Z}[\Lambda_{k}^{\ast}],
\]
where under the identification $\mathfrak{g}\equiv\mathfrak{g}^{\ast}$
given by the basic inner product, the level $k$ weights are given
by $\Lambda_{k}^{\ast}=\Lambda_{+}^{\ast}\cap k\mathfrak{A}$, and
the fusion ideal is defined as:
\[
I_{k}(G):=\left\{ f\in R(G)\mbox{ }\big|\mbox{ }f\left(\exp\left(\tfrac{\lambda+\rho}{k+\mathsf{h}^{\vee}}\right)\right)=0\mbox{ for all }\lambda\in\Lambda_{k}^{\ast}\right\} ,
\]
with $\mathsf{h}^{\vee}=1+\rho\cdot\theta$ the dual coxeter number
of $G$. We recall that for $G$ simply connected, $R_{k}(G)$ coincides
with the ring of level $k$ projective representations of $LG$, the
loop group associated to $G$ \cite[Rk. A.6]{[Mein12]}.

The main objects of this section are the $\kappa$-twisted analogues
of the rings $R(G)$ and $R_{k}(G)$ for $G$ simply connected and
simple, and we start by discussing the analogues of the irreducible
characters $\{\chi_{\lambda}\}_{\lambda\in\Lambda_{+}^{\ast}}$.

\subsection{Construction and First Properties}

For a dominant weight $\lambda\in\Lambda_{+}^{\ast}$, let $(\rho_{\lambda},V_{\lambda})$
denote the corresponding irreducible highest weight unitary representation,
with $\rho_{\lambda}:G\rightarrow\mathrm{U}(V_{\lambda})$ and normalized
highest weight vector $v_{\lambda}\in V_{\lambda}$, and let $(\Lambda_{+}^{\ast})^{\kappa}=(\Lambda^{\ast})^{\kappa}\cap(\mathfrak{t}^{\kappa})_{+}^{\ast}$
denote the $\kappa$-fixed integral dominant weights of $(G,T)$.
\begin{defn}
Let $(\rho,V)$ be a finite-dimensional unitary representation of
$G$. We call $(\rho,V)$ a $\kappa$\textbf{-admissible} representation
if it admits an implementation of $\kappa$ that is compatible with
the $G$-action, that is, if there exists a unitary operator $\tilde{\kappa}_{V}\in\mathrm{Aut}(V)$
such that:
\[
\rho_{V}\left(\kappa(g)\right)=\tilde{\kappa}_{V}\circ\rho_{V}(g)\circ\tilde{\kappa}_{V}^{-1},\ \ \forall g\in G.
\]
\end{defn}

We note that $\kappa$-admissible representations are closed under
direct sums and tensor products. For such a representation $(\rho,V)$,
the multiplicities of the weights $\{\kappa^{j}(\mu)\}_{j=1}^{|\kappa|}\subset\Lambda^{\ast}$
coincide, and its decomposition into irreducibles is of the form:
\begin{equation}
V=\bigoplus_{j=1}^{p_{0}}V_{\lambda_{j}^{(0)}}^{\oplus m_{j}^{(0)}}\oplus\bigoplus_{j=1}^{p_{1}}\left(\bigoplus_{l=1}^{|\kappa|}V_{\lambda_{j}^{(1)}}\right)^{\oplus m_{j}^{(1)}},\label{Eq_Kappa_Admissible_Repn_IrredDecomp}
\end{equation}
for some $m_{j}^{(k)},p_{k}\in\mathbb{N}$, $\{\lambda_{j}^{(0)}\}\subset(\Lambda_{+}^{\ast})^{\kappa}$
and $\{\lambda_{j}^{(1)}\}\subset\Lambda_{+}^{\ast}\smallsetminus(\Lambda_{+}^{\ast})^{\kappa}$.

For $\lambda\in(\Lambda_{+}^{\ast})^{\kappa}$, Schur's lemma gives
the existence of a unique unitary automorphism $\tilde{\kappa}_{\lambda}\in\mathrm{Aut}(V_{\lambda})$
such that:\begin{itemize}

\item[(i)] $\rho_{\lambda}\left(\kappa(g)\right)=\tilde{\kappa}_{\lambda}\circ\rho_{\lambda}(g)\circ\tilde{\kappa}_{\lambda}^{-1}$,
for all $g\in G$; 

\item[(ii)] $\tilde{\kappa}_{\lambda}(v_{\lambda})=v_{\lambda}$.

\end{itemize} For the weights $\lambda\in\Lambda_{+}^{\ast}\smallsetminus(\Lambda_{+}^{\ast})^{\kappa}$,
we use the fact that $(\rho_{\lambda}\circ\kappa,V_{\lambda})$ is
equivalent to $(\rho_{\kappa^{-1}(\lambda)},V_{\kappa^{-1}(\lambda)})$,
so that the automorphism $\kappa$ is implemented on $\oplus_{j}V_{\kappa^{j}(\lambda)}$
as a permutation of factors times a $|\kappa|$-th root of unity.
Following \cite{[FSS96],[Hon17a]}, we have:
\begin{defn}
Let $(\rho_{V},V)$ be a $\kappa$-admissible representation of $G$,
and let $\tilde{\kappa}_{V}\in\mathrm{Aut}(V)$ be an implementation
restricting to $\tilde{\kappa}_{\lambda}$ on each irreducible summand
$V_{\lambda}$ with $\lambda\in(\Lambda_{+}^{\ast})^{\kappa}$. The
$\bm{\kappa}$\textbf{-twining character} of $(\rho_{V},V)$ is the
function $G\rightarrow\mathbb{C}$ given for any $g\in G$ by:
\[
\tilde{\chi}_{V}^{\kappa}(g):=\mathrm{tr}_{_{V}}\left(\tilde{\kappa}_{V}\circ\rho_{V}(g)\right).
\]
\end{defn}

Note that for a $\kappa$-admissible representation $(\rho,V)$, the
non-zero contributions to the function $\tilde{\chi}_{V}^{\kappa}$
come only from the summands $V_{\lambda_{j}^{(0)}}$ in equation (\ref{Eq_Kappa_Admissible_Repn_IrredDecomp}),
where $\{\lambda_{j}^{(0)}\}\subset(\Lambda_{+}^{\ast})^{\kappa}$.
For a second admissible representation $(\rho',V')$, it is easily
verified that:
\[
\tilde{\chi}_{V\oplus V'}^{\kappa}=\tilde{\chi}_{V}^{\kappa}+\tilde{\chi}_{V'}^{\kappa},\ \ \tilde{\chi}_{V\otimes V'}^{\kappa}=\tilde{\chi}_{V}^{\kappa}\cdot\tilde{\chi}_{V'}^{\kappa}.
\]
We will use these facts implicitly throughout this section, and we
will be most interested in the case of $V=V_{\lambda}$ with $\lambda\in(\Lambda_{+}^{\ast})^{\kappa}$.
\begin{rem}
An alternative way of thinking of $\kappa$-admissible representations
calls upon the non-connected group $G\rtimes\langle\kappa\rangle$
\cite[\S 2.4]{[Mo00]}. Identifying the connected component $G\kappa\subseteq G\rtimes\langle\kappa\rangle$
with the group $G$ equipped with $\kappa$-twisted conjugation, a
$\kappa$-admissible representation of $G$ is the restriction of
a representation of $G\rtimes\langle\kappa\rangle$. By \cite[Prop.2.7]{[Mo00]},
an irreducible representation $V$ of $G\rtimes\langle\kappa\rangle$
is determined by a highest weight $\lambda\in\Lambda_{+}^{\ast}$
and the implementation $\tilde{\kappa}_{V}\in\mathrm{Aut}(V)$ of
$\kappa$:\begin{itemize}

\item[(a)] If $\kappa(\lambda)\ne\lambda$, then $V=\oplus_{j=0}^{|\kappa|-1}V_{\kappa^{j}(\lambda)}$,
and $\tilde{\kappa}_{V}$ permutes the factors $V_{\kappa^{j}(\lambda)}$.

\item[(b)] If $\kappa(\lambda)=\lambda$, then $V=V_{\lambda}$,
and $\tilde{\kappa}_{V}=e^{\frac{2\pi i}{|\kappa|}j}\tilde{\kappa}_{\lambda}$
for $j=1,\cdots,|\kappa|$, giving rise to $|\kappa|$ distinct irreps
$\rho_{\lambda,i}:G\rtimes\langle\kappa\rangle\to\mathrm{Aut}(V_{\lambda})$.

\end{itemize} Similarly, the twining characters on $G$ correspond
to restrictions of characters on $G\rtimes\langle\kappa\rangle$.
For $V$ irreducible \cite[Prop.2.8]{[Mo00]}, $\chi_{V}|_{G\kappa}\equiv0$
if $V$ is as in (a) above, while $\chi_{V}|_{G\kappa}\equiv e^{\frac{2\pi i}{|\kappa|}j}\tilde{\chi}_{\lambda}^{\kappa}$
for $V$ as in (b).
\end{rem}

Let $L^{2}(G)$ be the space of $\mathbb{C}$-valued $L^{2}$-functions
on $G$ with respect to the normalized Haar measure $dg$, and denote
by $\mathcal{C}^{\infty}(G)$ the smooth $\mathbb{C}$-valued functions.
We use $L^{2}(G\kappa)^{G}$ and $\mathcal{C}^{\infty}(G\kappa)^{G}$
to denote the subspaces of $\mathrm{Ad}^{\kappa}$-invariant functions
on $G$. As it would be expected, the irreducible twining characters
$\{\tilde{\chi}_{\lambda}^{\kappa}\}_{\lambda\in(\Lambda_{+}^{\ast})^{\kappa}}$
generalize several properties of the characters $\{\chi_{\lambda}\}_{\lambda\in\Lambda_{+}^{\ast}}$.
First, we remind the following facts on matrix coefficients \cite[\S\S  II.4,  III.3]{[BtD95]}:
\begin{lem}
\label{Lemma_Matrix_Coeff} For any finite-dimensional unitary representation
$(\rho_{V},V)$ of $G$, denote the matrix coefficient corresponding
to $\alpha\in V^{\ast}$ and $v\in V$ by: 
\[
m_{\alpha,v}^{V}:G\longrightarrow\mathbb{C}\mbox{, }g\longmapsto\alpha\left(\rho_{V}(g)\cdot v\right),
\]
and let $G$ act on $\mathrm{End}(V)$ via the assignment:
\[
g\cdot A:=\rho_{V}(g)\circ A\circ\rho_{V}(g)^{-1}\mbox{, for all }A\in\mathrm{End}(V).
\]
Denote by $M(G)$ the space of complex matrix coefficients on $G$.

\begin{enumerate}

\item The space $M(G)$ is dense in $\mathcal{C}^{0}(G)$ and in
$L^{2}(G)$.

\item For $(\rho_{V},V)$ irreducible, one has the following identity,
for all $A\in\mathrm{End}(V)$:
\[
\int_{G}dg\left(g\cdot A\right)=\frac{\mathrm{tr}_{_{V}}(A)}{\dim V}\mathrm{Id}_{V}.
\]

\item For irreducible representations $(\rho_{V},V)$ and $(\rho_{W},W)$,
the matrix coefficients satisfy the following orthogonality relations:
\[
\int_{G}\overline{\langle a,\rho_{V}(g)\cdot v\rangle_{_{V}}}\langle b,\rho_{W}(g)\cdot w\rangle_{_{W}}=\begin{cases}
\frac{1}{\dim V}\overline{\langle a,b\rangle}_{_{V}}\langle v,w\rangle_{_{V}}, & \mbox{if }V\simeq W;\\
0, & \mbox{if }V\mbox{\ensuremath{\not\simeq}}W,
\end{cases}
\]
for any $a,v\in V$ and $b,w\in W$, and $\langle-,-\rangle_{_{V}}$
denoting the Hermitian inner product on $V$.

\end{enumerate}
\end{lem}

The main result of this subsection is the following:
\begin{prop}
\label{Prop_Twining_Char_L2_Class_Fns}With the notations of this
section: \begin{enumerate}

\item The averaging map $\mathrm{Av}^{\kappa}:L^{2}(G)\rightarrow L^{2}(G\kappa)^{G}\mbox{, }f\mapsto\int_{G}dg\cdot(\mathrm{Ad}_{g}^{\kappa})^{\ast}f$
is an orthogonal projection.

\item The twining characters $\{\tilde{\chi}_{\lambda}^{\kappa}\}_{\lambda\in(\Lambda_{+}^{\ast})^{\kappa}}$
generate a dense subspace of $L^{2}(G\kappa)^{G}$, and satisfy the
orthogonality relations:
\[
\langle\tilde{\chi}_{\lambda}^{\kappa},\tilde{\chi}_{\mu}^{\kappa}\rangle_{L^{2}}=\delta_{\mu\lambda},\ \ \lambda,\mu\in(\Lambda_{+}^{\ast})^{\kappa}.
\]

\item For a second Dynkin diagram automorphism $\tau\in\mathrm{Out}(G)$,
one has: 
\[
\left(\tilde{\chi}_{\lambda}^{\kappa}\ast\tilde{\chi}_{\mu}^{\tau}\right)(x)=(\dim V_{\lambda})^{-1}\delta_{\lambda\mu}\tilde{\chi}_{\lambda}^{\tau\kappa}(x)=(\dim V_{\lambda})^{-1}\delta_{\lambda\mu}\tilde{\chi}_{\lambda}^{\kappa\tau}\left(\kappa(x)\right),
\]
for any $\lambda,\mu\in(\Lambda_{+}^{\ast})^{\kappa}\cap(\Lambda_{+}^{\ast})^{\tau}$,
where convolution is given by: 
\[
\left(\psi\ast\varphi\right)(x)=\int_{G}dg\cdot\psi(xg^{-1})\varphi(g),\ \ \psi,\varphi\in L^{2}(G).
\]

\item For $\lambda\in(\Lambda_{+}^{\ast})^{\kappa}$, the ``spherical
harmonic'' $\phi_{\lambda}:G\rightarrow\mathbb{C}$, $g\mapsto\langle v_{\lambda},\rho_{\lambda}(g)\cdot v_{\lambda}\rangle_{V_{\lambda}}$
satisfies the identity: 
\[
\int_{G}dg\cdot\phi_{\lambda}\left(\mathrm{Ad}_{g}^{\kappa}x\right)=(\dim V_{\lambda})^{-1}\tilde{\chi}_{\lambda}^{\kappa}(x).
\]

\end{enumerate}
\end{prop}

\begin{proof}
\noindent \begin{enumerate}

\item The formula in the statement defines a bounded linear operator
$\mathrm{Av}^{\kappa}:L^{2}(G)\rightarrow L^{2}(G)$ with image in
$L^{2}(G\kappa)^{G}$. That it is an orthogonal projection follows
from the fact that it is self-adjoint and coincides with the identity
on $L^{2}(G\kappa)^{G}$.

\item That the $\kappa$-twining characters $\{\tilde{\chi}_{\lambda}^{\kappa}\}_{\lambda\in(\Lambda_{+}^{\ast})^{\kappa}}$
are $\mathrm{Ad}^{\kappa}$-invariant follows from the cyclic property
of the trace and the defining identity:
\[
\rho_{\lambda}\left(\kappa(g)\right)=\tilde{\kappa}_{\lambda}\circ\rho_{\lambda}(g)\circ\tilde{\kappa}_{\lambda}^{-1}.
\]
Next, let $\alpha\in V_{\lambda}^{\ast}$, $v\in V_{\lambda}$ and
$A=\tilde{\kappa}_{\lambda}\circ\rho_{\lambda}(g)$. Using lemma \ref{Lemma_Matrix_Coeff},
we have for any $x\in G$ that:
\begin{eqnarray}
(\mathrm{Av}^{\kappa}m_{\alpha,v}^{V_{\lambda}})(x) & = & \alpha\left(\int_{G}\left(\rho_{\lambda}(\mathrm{Ad}_{g}^{\kappa}x)\cdot v\right)dg\right)=\alpha\left(\int_{G}(g\cdot A)\left(\tilde{\kappa}_{\lambda}^{-1}\cdot v\right)\right)\nonumber \\
 & = & \alpha\left(\frac{\mathrm{tr}_{_{V_{\lambda}}}(A)}{\dim V_{\lambda}}\tilde{\kappa}_{\lambda}^{-1}\cdot v\right)=\frac{\alpha\left(\tilde{\kappa}_{\lambda}^{-1}\cdot v\right)}{\dim V_{\lambda}}\tilde{\chi}_{\lambda}^{\kappa}(x).\label{Eq_Avg_Elem_Matrix_Coeff}
\end{eqnarray}
This is used for the density of twining characters as follows. Given
$\varphi\in L^{2}(G\kappa)^{G}$, for any $\varepsilon>0$ there is
a function $f\in M(G)$ such that $|\!|\varphi-f|\!|_{L^{2}}<\varepsilon$,
and such that for a finite subset $L\subseteq\Lambda_{+}^{\ast}$:
\[
f(x)=\sum_{\lambda\in L}\sum_{i=1}^{k_{\lambda}}a_{i}m_{\alpha_{i},v_{i}}^{V_{\lambda}}(x),
\]
where $\{k_{\lambda}\}_{\lambda\in L}$ are integers, and for all
$1\le i\le k_{\lambda}$: $a_{i}\in\mathbb{C}$, $v_{i}\in V_{\lambda}$,
$\alpha_{i}\in V_{\lambda}^{\ast}$. By equation (\ref{Eq_Avg_Elem_Matrix_Coeff}),
we have:
\[
\mathrm{Av}^{\kappa}f(x)=\sum_{\lambda\in L\cap(\Lambda_{+}^{\ast})^{\kappa}}c_{\lambda}\tilde{\chi}_{\lambda}^{\kappa}(x),
\]
for some constants $c_{\lambda}\in\mathbb{C}$, and our initial function
$\varphi\in L^{2}(G\kappa)^{G}$ then satisfies the estimate:
\[
|\!|\varphi-\mathrm{Av}^{\kappa}f|\!|_{L^{2}}=|\!|\mathrm{Av}^{\kappa}(\varphi-f)|\!|_{L^{2}}\le|\!|\varphi-f|\!|_{L^{2}}<\varepsilon,
\]
which establishes the density claim. For the orthogonality relations,
let $\{v_{i}\}\subset V_{\lambda}$ and $\{w_{i}\}\subset V_{\mu}$
be orthonormal bases, so that:
\[
\langle\tilde{\chi}_{\lambda}^{\kappa},\tilde{\chi}_{\mu}^{\kappa}\rangle_{L^{2}}=\sum_{i,j}\int_{G}\overline{\langle v_{i},(\tilde{\kappa}_{\lambda}g)\cdot v_{i}\rangle}_{V_{\lambda}}\langle w_{j},(\tilde{\kappa}_{\mu}g)\cdot w_{j}\rangle_{_{V_{\mu}}}.
\]
 By lemma \ref{Lemma_Matrix_Coeff}-(3), the above vanishes if $\mu\ne\lambda$,
and otherwise:
\[
\langle\tilde{\chi}_{\lambda}^{\kappa},\tilde{\chi}_{\lambda}^{\kappa}\rangle_{L^{2}}=\frac{1}{\dim V_{\lambda}}\sum_{i,j}\overline{\langle v_{i},v_{j}\rangle}_{V_{\lambda}}\langle(\tilde{\kappa}_{\lambda}g)\cdot v_{i},(\tilde{\kappa}_{\lambda}g)\cdot v_{j}\rangle_{_{V_{\lambda}}}=1.
\]

\item Let $\lambda,\mu\in(\Lambda_{+}^{\ast})^{\kappa}\cap(\Lambda_{+}^{\ast})^{\tau}$
and let $\{v_{i}\}\subset V_{\lambda}$ and $\{w_{i}\}\subset V_{\mu}$
be orthonormal bases. If $\mu\ne\lambda$, the integral: 
\[
(\tilde{\chi}_{\lambda}^{\kappa}\ast\tilde{\chi}_{\mu}^{\tau})(x)=\sum_{i,j}\int_{G}\langle v_{i},(\tilde{\kappa}_{\lambda}xg^{-1})\cdot v_{i}\rangle_{_{V_{\lambda}}}\langle w_{i},(\tilde{\tau}_{\mu}xg^{-1})\cdot w_{i}\rangle_{_{V_{\mu}}}dg
\]
vanishes by \ref{Lemma_Matrix_Coeff}-(3). To compute $\tilde{\chi}_{\lambda}^{\kappa}\ast\tilde{\chi}_{\lambda}^{\tau}$,
we use $\tilde{\chi}_{\lambda}^{\kappa}(x)=\mathrm{tr}_{_{V_{\lambda}}}(x\tilde{\kappa}_{\lambda})$
and notice that: 
\[
\langle v_{i},(xg^{-1}\tilde{\kappa}_{\lambda})\cdot v_{i}\rangle_{_{V_{\lambda}}}=\sum_{k=1}^{\dim V_{\lambda}}\langle v_{i},x\cdot v_{k}\rangle_{_{V_{\lambda}}}\overline{\langle\tilde{\kappa}_{\lambda}\cdot v_{i},g\cdot v_{k}\rangle}_{V_{\lambda}}.
\]
This identity along with lemma \ref{Lemma_Matrix_Coeff}-(3) yield:
\begin{eqnarray*}
(\tilde{\chi}_{\lambda}^{\kappa}\ast\tilde{\chi}_{\lambda}^{\tau})(x) & = & \sum_{i,j}\int_{G}\langle v_{i},(xg^{-1}\tilde{\kappa}_{\lambda})\cdot v_{i}\rangle_{_{V_{\lambda}}}\langle v_{j},(g\tilde{\tau}_{\lambda})v_{j}\rangle_{_{V_{\lambda}}}dg\\
 & = & \sum_{i,j,k}\langle v_{i},x\cdot v_{k}\rangle_{_{V_{\lambda}}}\int_{G}\overline{\langle\tilde{\kappa}_{\lambda}\cdot v_{i},g\cdot v_{k}\rangle}_{V_{\lambda}}\langle v_{j},(g\tilde{\tau}_{\lambda})v_{j}\rangle_{_{V_{\lambda}}}dg\\
 & = & \frac{1}{\dim V_{\lambda}}\sum_{i,j,k}\langle v_{i},x\cdot v_{k}\rangle_{_{V_{\lambda}}}\langle v_{j},\tilde{\kappa}_{\lambda}\cdot v_{i}\rangle_{_{V_{\lambda}}}\langle v_{k},\tilde{\tau}_{\lambda}\cdot v_{j}\rangle_{_{V_{\lambda}}}\\
 & = & \frac{1}{\dim V_{\lambda}}\sum_{i}\langle v_{i},(x\tilde{\tau}_{\lambda}\tilde{\kappa}_{\lambda})\cdot v_{i}\rangle_{_{V_{\lambda}}}=\frac{1}{\dim V_{\lambda}}\mathrm{tr}_{V_{\lambda}}\left(\rho_{\lambda}(x)\circ(\tilde{\tau}_{\lambda}\tilde{\kappa}_{\lambda})\right).
\end{eqnarray*}
The operator $(\tilde{\tau}_{\lambda}\tilde{\kappa}_{\lambda})\in\mathrm{Aut}(V_{\lambda})$
preserves the highest weight vector of $V_{\lambda}$ and satisfies:
\[
\rho_{\lambda}\left(\tau\kappa(g)\right)=(\tilde{\tau}_{\lambda}\tilde{\kappa}_{\lambda})\circ\rho_{\lambda}(g)\circ(\tilde{\tau}_{\lambda}\tilde{\kappa}_{\lambda})^{-1},\ \ g\in G.
\]
By uniqueness, we must have $(\tilde{\tau}_{\lambda}\tilde{\kappa}_{\lambda})=\widetilde{(\tau\kappa)}_{\lambda}$,
and therefore: 
\[
\tilde{\chi}_{\lambda}^{\kappa}\ast\tilde{\chi}_{\mu}^{\tau}=(\dim V_{\lambda})^{-1}\delta_{\lambda\mu}\tilde{\chi}_{\lambda}^{\tau\kappa}.
\]
That one has $\tilde{\chi}_{\lambda}^{\tau\kappa}(g)=\tilde{\chi}_{\lambda}^{\kappa\tau}\left(\kappa(g)\right)$
for all $g\in G$ and $\lambda\in(\Lambda_{+}^{\ast})^{\kappa}\cap(\Lambda_{+}^{\ast})^{\tau}$
follows from:
\[
\mathrm{tr}_{_{V_{\lambda}}}\left((\tilde{\tau}_{\lambda}\tilde{\kappa}_{\lambda})\rho_{\lambda}(g)\right)=\mathrm{tr}_{_{V_{\lambda}}}\left(\tilde{\tau}_{\lambda}\rho_{\lambda}\left(\kappa(g)\right)\tilde{\kappa}_{\lambda}\right)=\mathrm{tr}_{_{V_{\lambda}}}\left((\tilde{\kappa}_{\lambda}\tilde{\tau}_{\lambda})\rho_{\lambda}\left(\kappa(g)\right)\right).
\]

\item The identity of the statement is a special case of equation
(\ref{Eq_Avg_Elem_Matrix_Coeff}) with $v=v_{\lambda}=\tilde{\kappa}_{\lambda}(v_{\lambda})$
and $\alpha=\langle v_{\lambda},-\rangle_{_{V_{\lambda}}}$.\qedhere

\end{enumerate}
\end{proof}

\subsection{The Jantzen Character Formula}

We now come to a Weyl-type character formula for the $\widetilde{\chi}_{\lambda}^{\kappa}$,
for which we will need an integration formula. If $(T^{\kappa})^{\kappa-\mathrm{reg}}$
and $G^{\kappa-\mathrm{reg}}$ denote the submanifolds of $\mathrm{Ad}^{\kappa}$-regular
elements in $T^{\kappa}$ and $G$ (elements $g$ with $\dim Z_{G}^{\kappa}(g)=\dim T^{\kappa}$),
then proposition \ref{Prop_Twisted_Conj_Classes} and lemma \ref{Lemma_Det_Diff_Conj_Map}
imply that the restriction of the conjugation map $\mathsf{c}$ to
$(T^{\kappa})^{\kappa-\mathrm{reg}}\times(G/T^{\kappa})$ is a covering
onto $G^{\kappa-\mathrm{reg}}$, with group of deck transformations
$W^{(\kappa)}$. Using Fubini's theorem along with the fact that the
connected components of $G\smallsetminus G^{\kappa-\mathrm{reg}}$
and $T^{\kappa}\smallsetminus(T^{\kappa})^{\kappa-\mathrm{reg}}$
have measure zero, we obtain \cite[Thm.2.5]{[Wen01]}:
\begin{lem}
\label{Lemma_Twisted_Integral_Formula}For a class function $f\in L^{1}(G\kappa)^{G}$,
one has that:
\begin{eqnarray*}
\int_{G}f(x)dx & = & \frac{1}{|W^{\kappa}|}\int_{T^{\kappa}}f(t)\big|\widetilde{\Delta}(t)\big|^{2}dt.
\end{eqnarray*}
\end{lem}

\noindent We will also need some observations on the action of the
induced automorphism $\tilde{\kappa}_{\lambda}\in\mathrm{Aut}(V_{\lambda})$
on the highest weight irrep $V_{\lambda}$, which are easily established
using the weight space decomposition:
\begin{lem}
\label{Lemma_Action_kappa_Weights_Irreps}For $\lambda\in(\Lambda_{+}^{\ast})^{\kappa}$,
let $P(V_{\lambda})$ be the set of weights of the representation
$V_{\lambda}$, and let $V_{\lambda}^{\kappa}=\ker(\tilde{\kappa}_{\lambda}-1)$
and $P(V_{\lambda})^{\kappa}=\{\mu\in P(V_{\lambda})\mbox{ }|\mbox{ }\kappa(\mu)=\mu\}$.

\begin{enumerate}

\item The set $P(V_{\lambda})\subset\Lambda^{\ast}$ is preserved
by $\kappa\in\mathrm{Out}(G)$, and if $u\in V_{\lambda}$ is a weight
vector for $\mu\in P(V_{\lambda})$, then $\tilde{\kappa}_{\lambda}(u)\in V_{\lambda}$
is a weight vector for $\kappa(\mu)\in P(V_{\lambda})$.

\item For $t\in T^{\kappa}$, one has that:
\[
\widetilde{\chi}_{\lambda}^{\kappa}(t)=\mathrm{tr}_{V_{\lambda}^{\kappa}}\left(\rho_{\lambda}(t)\right)=\sum_{\mu\in P(V_{\lambda})^{\kappa}}m_{\mu}t^{\mu},
\]
with $m_{\mu}=\dim(V_{\lambda})_{\mu}$, and that: 
\[
\widetilde{\chi}_{\lambda}^{\kappa}(e)=\dim V_{\lambda}^{\kappa}.
\]

\end{enumerate}
\end{lem}

For a weight $\lambda\in(\Lambda_{+}^{\ast})^{\kappa}$, let $\tilde{V}_{\lambda}$
denote the corresponding highest weight irrep of the orbit group $G_{(\kappa)}$,
and let $\sigma_{\lambda}:G_{(\kappa)}\rightarrow\mathbb{C}$ be the
associated irreducible character. Recalling that the half-sum of positive
roots of $\mathfrak{R}_{(\kappa)}$ coincides with $\rho=\frac{1}{2}\sum_{\alpha\in\mathfrak{R}_{+}}\alpha$,
we consider the functions:
\[
\widetilde{J}_{\lambda}(\xi):=\sum_{w\in W^{\kappa}}(-1)^{l(w)}e^{2\pi i\langle w\cdot(\lambda+\rho),\xi\rangle},\ \ \xi\in\mathfrak{t}^{\kappa},
\]
so that by the Weyl character formula, for any $\xi\in\mathfrak{t}^{\kappa}$
regular:
\[
\sigma_{\lambda}\left(\exp_{\mathfrak{g}_{(\kappa)}}(\xi)\right)=\frac{\widetilde{J}_{\lambda}(\xi)}{\widetilde{J}_{0}(\xi)},
\]
where $\exp_{\mathfrak{g}_{(\kappa)}}:\mathfrak{g}_{(\kappa)}\rightarrow G_{(\kappa)}$
sends $\mathfrak{t}^{\kappa}$ to $T_{(\kappa)}=T^{\kappa}/(T^{\kappa}\cap T_{\kappa})$
(we make this distinction since $\exp_{\mathfrak{g}}:\mathfrak{g}\rightarrow G$
is also used below, and sends $\mathfrak{t}^{\kappa}$ to $T^{\kappa}$).
The next result was first proved by Jantzen in \cite{[Jant73]}, then
recovered in the context of affine algebras in \cite[Thm.1]{[FSS96]},
and in the context of non-connected groups in \cite[Thm.2.6]{[Wen01]}.
\begin{prop}
\textbf{\textup{Jantzen Character Formula\label{Prop_Jantzen_Char_Formula}}}
Let $G$ be a compact, connected, simply connected and simple Lie
group, let $\kappa\in\mathrm{Out}(G)$ be a Dynkin diagram automorphism
and $G_{(\kappa)}$ the associated orbit Lie group. For an invariant
dominant weight $\lambda\in(\Lambda_{+}^{\ast})^{\kappa}$, one has
that:
\[
\left(\widetilde{\chi}_{\lambda}^{\kappa}\right)_{|T^{\kappa}}=\pi^{\ast}\left(\sigma_{\lambda|T_{(\kappa)}}\right),
\]
where $\pi:T^{\kappa}\rightarrow T_{(\kappa)}$ is the quotient map
by the action of $T^{\kappa}\cap T_{\kappa}$. In particular, for
$\xi\in\mathfrak{t}^{\kappa}$ \textup{regular}:
\[
\widetilde{\chi}_{\lambda}^{\kappa}\left(\exp_{\mathfrak{g}}(\xi)\right)=\frac{\sum_{w\in W^{\kappa}}(-1)^{l(w)}e^{2\pi i\langle w\cdot(\lambda+\rho)-\rho,\xi\rangle}}{\widetilde{\Delta}\left(\exp_{\mathfrak{g}}(\xi)\right)}=\frac{\widetilde{J}_{\lambda}(\xi)}{\widetilde{J}_{0}(\xi)}.
\]
\end{prop}

\begin{proof}
Since $\widetilde{\chi}_{\lambda}^{\kappa}$ is $\mathrm{Ad}^{\kappa}$-invariant,
$\left(\widetilde{\chi}_{\lambda}^{\kappa}\right)_{|T^{\kappa}}$
is invariant under the action of $T^{\kappa}\cap T_{\kappa}$ by multiplication.
Along with lemma \ref{Lemma_Action_kappa_Weights_Irreps}-(2), we
thus have for any $t\in T^{\kappa}$ that:
\[
\left(\widetilde{\chi}_{\lambda}^{\kappa}\right)_{|T^{\kappa}}(t)=\sum_{\mu\in P(V_{\lambda})^{\kappa}}m_{\mu}t^{\mu}=\sum_{\mu\in P(V_{\lambda})^{\kappa}}m_{\mu}\pi(t)^{\mu}.
\]
On the one hand, this shows that $P(V_{\lambda})^{\kappa}$ is a subset
of the weight lattice $(\Lambda^{\ast})^{\kappa}$ of $G_{(\kappa)}$,
and on the other hand, it means that $\left(\widetilde{\chi}_{\lambda}^{\kappa}\right)_{|T^{\kappa}}=\pi^{\ast}(\sigma_{\tilde{V}|T_{(\kappa)}})$,
where $\sigma_{\tilde{V}}:G_{(\kappa)}\rightarrow\mathbb{C}$ is the
character of a $G_{(\kappa)}$-module $\tilde{V}=\bigoplus_{\mu\in(\Lambda_{+}^{\ast})^{\kappa}}\tilde{V}_{\mu}^{\oplus q_{\mu}}$:
\[
\sigma_{\tilde{V}}=\sum_{\mu\in(\Lambda_{+}^{\ast})^{\kappa}}q_{\mu}\sigma_{\mu}\mbox{, }
\]
with finitely many nonzero integers $q_{\mu}$. Since the function
$\widetilde{\Delta}(t)=\prod_{\alpha\in\mathfrak{R}_{(\kappa)+}}(1-t^{-\alpha})$
is the ``Weyl denominator'' on $T_{(\kappa)}$, we have by lemma
\ref{Lemma_Twisted_Integral_Formula} and the orthogonality relations
for twining characters that:
\[
\sum_{\mu\in(\Lambda_{+}^{\ast})^{\kappa}}q_{\mu}^{2}=\frac{1}{|W^{\kappa}|}\int_{T_{(\kappa)}}ds\cdot\big|\widetilde{\Delta}\sigma_{\tilde{V}}\big|^{2}=\frac{1}{|W^{\kappa}|}\int_{T^{\kappa}}dt\cdot\big|\widetilde{\Delta}\left(\widetilde{\chi}_{\lambda}^{\kappa}\right)_{|T^{\kappa}}\big|^{2}=1.
\]
We thus have $\tilde{V}=\tilde{V}_{\lambda}\simeq V_{\lambda}^{\kappa}$,
since the equation above shows that $\tilde{V}$ is irreducible, and
since $\lambda\in(\Lambda_{+}^{\ast})^{\kappa}$ is the highest weight
in $P(V_{\lambda})^{\kappa}$. This shows that $\left(\widetilde{\chi}_{\lambda}^{\kappa}\right)_{|T^{\kappa}}=\pi^{\ast}\left(\sigma_{\lambda|T_{(\kappa)}}\right)$,
and for any $\xi\in\mathfrak{t}^{\kappa}$ regular, we have by the
character formula for $\sigma_{\lambda}$ that: 
\begin{eqnarray*}
\sigma_{\lambda}\left(\exp_{\mathfrak{g}_{(\kappa)}}(\xi)\right) & = & \frac{\widetilde{J}_{\lambda}(\xi)}{\widetilde{J}_{0}(\xi)}=\frac{\sum_{w\in W^{\kappa}}(-1)^{l(w)}e^{2\pi i\langle w\cdot(\lambda+\rho),\xi\rangle}}{\sum_{w\in W^{\kappa}}(-1)^{l(w)}e^{2\pi i\langle w\cdot\rho,\xi\rangle}}\\
 & = & \frac{\sum_{w\in W^{\kappa}}(-1)^{l(w)}e^{2\pi i\langle w\cdot(\lambda+\rho)-\rho,\xi\rangle}}{\sum_{w\in W^{\kappa}}(-1)^{l(w)}e^{2\pi i\langle w\cdot\rho-\rho,\xi\rangle}}\\
 & = & \frac{\sum_{w\in W^{\kappa}}(-1)^{l(w)}e^{2\pi i\langle w\cdot(\lambda+\rho)-\rho,\xi\rangle}}{\widetilde{\Delta}\left(\exp_{\mathfrak{g}}(\xi)\right)},
\end{eqnarray*}
which is the second equation of the proposition.
\end{proof}
The Jantzen character formula shows that $\pi^{\ast}:L^{2}(T_{(\kappa)})^{W^{\kappa}}\rightarrow L^{2}(T^{\kappa})^{W^{(\kappa)}}$
is an isomorphism, and this is true for smooth class functions in
particular. With a little more work, we also have the following version
of the Chevalley restriction theorem:
\begin{prop}
Let $G$ be a compact, connected, simply connected and simple Lie
group, and $\kappa\in\mathrm{Out}(G)$ a Dynkin diagram automorphism.
The inclusion $\iota_{T^{\kappa}}:T^{\kappa}\hookrightarrow G$ induces
the following isomorphism of algebras:
\[
\iota_{T^{\kappa}}^{\ast}:\mathcal{C}^{\infty}(G\kappa)^{G}\longrightarrow\mathcal{C}^{\infty}(T^{\kappa})^{W^{(\kappa)}}\simeq\mathcal{C}^{\infty}(T_{(\kappa)})^{W^{\kappa}}\mbox{, }f\longmapsto f_{|T^{\kappa}}
\]
where $\mathcal{C}^{\infty}(G\kappa)^{G}$ is the ring of smooth $\mathrm{Ad}^{\kappa}$-invariant
functions on $G$ with values in $\mathbb{C}$.
\end{prop}

\begin{proof}
Using a $\mathrm{Ad}^{\kappa}(G)$-invariant partition of unity on
$G$, the problem is reduced to proving that for any $x\in G$, there
exists a $\mathrm{Ad}^{\kappa}(G)$-invariant open neighborhood $\mathcal{U}_{x}$
such that the restriction map: 
\[
\mathcal{C}_{c}^{\infty}(\mathcal{U}_{x}\kappa)^{G}\rightarrow\mathcal{C}_{c}^{\infty}(\mathcal{U}_{x}\cap T^{\kappa})^{W^{(\kappa)}},\ f\mapsto f_{|\mathcal{U}_{x}\cap T^{\kappa}},
\]
is an isomorphism. To obtain the open subset $\mathcal{U}_{x}$, we
construct a slice for the twisted adjoint action following the approach
of \cite[Prop.2.5]{[Mein17]}, and to lighten the notation, we write
$Z_{x}^{\kappa}$ for the stabilizer of $x$ with Lie algebra $\mathfrak{z}_{x}^{\kappa}$.

Let $F\subset\mathfrak{A}_{(\kappa)}$ be an open face of the fundamental
alcove in $\mathfrak{t}^{\kappa}$, and let $\xi\in F$ with $x=e^{\xi}$.
Let $V_{\xi}\subset\mathfrak{t}^{\kappa}$ be an open ball centred
at $\xi$, and let $\mathring{\mathfrak{A}}_{(\kappa)}^{F}$ be the
union of open faces in $\mathfrak{A}_{(\kappa)}$ containing $F$
in their closure. We define the following $\mathrm{Ad}^{\kappa}\left(Z_{x}^{\kappa}\right)$-invariant
open subset of $Z_{x}^{\kappa}$: 
\[
U_{x}:=\mathrm{Ad}^{\kappa}\left(Z_{x}^{\kappa}\right)\cdot\exp\left(V_{\xi}\cap\mathring{\mathfrak{A}}_{(\kappa)}^{F}\right),
\]
which for $V_{\xi}$ small enough gives a $Z_{x}^{\kappa}$-equivariant
diffeomorphism :
\[
x\exp_{\mathfrak{z}_{x}^{\kappa}}:V_{0}\longrightarrow U_{x}\mbox{, }\zeta\longmapsto xe^{\zeta}=e^{\kappa(\zeta)}x,
\]
with $V_{0}\subset\mathfrak{z}_{x}^{\kappa}$ a $Z_{x}^{\kappa}$-invariant
open ball centred at the origin, and such that:
\[
x\exp_{\mathfrak{z}_{x}^{\kappa}}\left(\mathrm{Ad}_{\kappa(g)}\zeta\right)=\mathrm{Ad}_{g}^{\kappa}\left(x\exp_{\mathfrak{z}_{x}^{\kappa}}(\zeta)\right),\ \ \forall g\in Z_{x}^{\kappa}.
\]
(In the left trivialization of $TG$, the twisted adjoint action of
$g\in G$ on $\zeta\in\mathfrak{g}$ is given by $\mathrm{Ad}_{\kappa(g)}\zeta$.)
By construction, $U_{x}\subset Z_{x}^{\kappa}$ is a slice through
$x$ for $\kappa$-twisted conjugation, and the map:
\[
G\times_{Z_{x}^{\kappa}}V_{0}\longrightarrow G\mbox{, }[g,\zeta]\longmapsto\mathrm{Ad}_{g}^{\kappa}\left(xe^{\zeta}\right)
\]
gives a $\mathrm{Ad}^{\kappa}(G)$-equivariant diffeomorphism onto
its image, which we denote henceforth by $\mathcal{U}_{x}$. Recalling
that $T^{\kappa}$ is a maximal torus of $Z_{x}^{\kappa}$, we also
obtain a $W_{x}^{(\kappa)}=N_{Z_{x}^{\kappa}}^{\kappa}(T^{\kappa})/T^{\kappa}$-equivariant
diffeomorphism:

\[
N_{G}^{\kappa}(T^{\kappa})\times_{N_{Z_{x}^{\kappa}}^{\kappa}(T^{\kappa})}\left(V_{0}\cap\mathfrak{t}^{\kappa}\right)\longrightarrow U_{x}\cap T^{\kappa}\mbox{, }[g,\zeta]\longmapsto\mathrm{Ad}_{g}^{\kappa}\left(xe^{\zeta}\right).
\]
By the Chevalley restriction theorem for Lie algebras \cite[Thm.7.28]{[BGV]},
we have an isomorphism $\mathcal{C}^{\infty}(\mathfrak{z}_{x}^{\kappa})^{Z_{x}^{\kappa}}\simeq\mathcal{C}^{\infty}(\mathfrak{t}^{\kappa})^{W_{x}^{(\kappa)}}$,
and therefore 
\[
\mathcal{C}_{c}^{\infty}(V_{0})^{Z_{x}^{\kappa}}\simeq\mathcal{C}_{c}^{\infty}(V_{0}\cap\mathfrak{t}^{\kappa})^{W_{x}^{(\kappa)}},
\]
and the equivariant diffeomorphisms constructed above yield:
\[
\mathcal{C}_{c}^{\infty}(\mathcal{U}_{x}\kappa)^{G}\simeq\mathcal{C}_{c}^{\infty}(\mathcal{U}_{x}\cap T^{\kappa})^{W^{(\kappa)}}.\qedhere
\]
\end{proof}

\subsection{Twisted Representation and Fusion Rings}

We turn to the representation rings associated to the $\kappa$-admissible
representations of $G$ in this subsection, following the approach
of \cite[App.4]{[Mein12]} (proofs can be found in \cite[Pt.III]{[Be93]}).
\begin{defn}
Define the \textbf{twining representation ring} $\tilde{R}^{(\kappa)}(G)$
to be the subring of $\mathcal{C}^{\infty}(G\kappa)^{G}$ generated
by the twining characters of finite-dimensional $\kappa$-admissible
representations of $G$.
\end{defn}

\noindent This is the ring of virtual $\kappa$-admissible representations
of $G$. As a $\mathbb{Z}$-module, it is spanned by the irreducible
twining characters:
\[
\tilde{R}^{(\kappa)}(G)\simeq\mathbb{Z}\left[(\Lambda_{+}^{\ast})^{\kappa}\right].
\]
The last two results of the previous subsection yield the following
fact:
\begin{prop}
Let $G$ be a compact, connected, simply connected and simple Lie
group, $\kappa\in\mathrm{Out}(G)$ a diagram automorphism and $G_{(\kappa)}$
the corresponding orbit Lie group. One then has:
\[
\tilde{R}^{(\kappa)}(G)\simeq\tilde{R}^{(\kappa)}(T^{\kappa})^{W^{(\kappa)}}\simeq R\left(T_{(\kappa)}\right)^{W^{\kappa}}\simeq R\left(G_{(\kappa)}\right).
\]
\end{prop}

\noindent As in the case of the usual representation ring $R(G)$,
the twining representation ring is equipped with an involution:
\[
(\cdot)^{\ast}:\tilde{R}^{(\kappa)}(G)\longrightarrow\tilde{R}^{(\kappa)}(G)\mbox{, }
\]
which acts as $\tilde{\chi}_{V}^{(\kappa)}\longmapsto\tilde{\chi}_{V^{\ast}}^{(\kappa)}$
for any $\kappa$-admissible representation $V$ of $G$. There is
also a trace map:
\[
\mathrm{Tr}_{0}:\tilde{R}^{(\kappa)}(G)\longrightarrow\mathbb{Z}\mbox{, }f\longmapsto\int_{G}dg\cdot\left(\overline{\tilde{\chi}_{0}^{\kappa}}f\right),
\]
that gives the coefficient of the basis element corresponding to the
weight $\lambda=0$.

Now, consider the following rescaling of the inner product $B$ on
$\mathfrak{g}$:
\[
\tilde{B}:=\tfrac{2}{B\left(\theta_{(\kappa),\mathrm{l}}^{\vee},\theta_{(\kappa),\mathrm{l}}^{\vee}\right)}B,
\]
where $\theta_{(\kappa),\mathrm{l}}\in\mathfrak{R}_{(\kappa)+}$ is
the highest root of $G_{(\kappa)}$. The restriction $\tilde{B}_{|\mathfrak{t}^{\kappa}}\in S^{2}(\mathfrak{t}^{\kappa})^{\ast}$
coincides with the restriction of the \textit{basic inner-product}
on $\mathfrak{g}_{(\kappa)}$ to $\mathfrak{t}^{\kappa}$, which is
the only inner-product on $\mathfrak{t}^{\kappa}$ such that $\tilde{B}^{-1}(\tilde{\alpha},\tilde{\alpha})=2$
for any long root $\tilde{\alpha}\in\mathfrak{R}_{(\kappa)}$ and
restricting to a $\mathbb{Z}$-valued bilinear form on $\Lambda_{(\kappa)}$.
In particular, the isomorphism:
\[
\tilde{B}^{\flat}:\mathfrak{g}\longrightarrow\mathfrak{g}^{\ast}\mbox{, }\xi\longmapsto\tilde{B}(\xi,\cdot)
\]
maps $\Lambda_{(\kappa)}$ onto $(\Lambda^{\ast})^{\kappa}$, since
$\tilde{B}_{|\mathfrak{t}^{\kappa}}\in(\Lambda^{\ast})^{\kappa}\otimes(\Lambda^{\ast})^{\kappa}$.
Fixing a positive integer $k\in\mathbb{Z}_{>0}$, we introduce the
following objects:\begin{itemize}

\item The isomorphism $\tilde{B}_{k}^{\flat}:=k\tilde{B}^{\flat}$,
with inverse:

\[
\tilde{B}_{k}^{\sharp}:\mathfrak{g}{}^{\ast}\longrightarrow\mathfrak{g}\mbox{, }\lambda\longmapsto\tfrac{1}{k}(\tilde{B}^{\flat})^{-1}(\lambda).
\]

\item The $k$-rescaled alcove in $(\mathfrak{t}^{\kappa})^{\ast}$
(proposition \ref{Prop_Twisted_Alcove}):
\[
\mathfrak{A}_{(\kappa),k}^{\ast}:=\tilde{B}_{k}^{\flat}\left(\mathfrak{A}_{(\kappa)}\right).
\]

\item The set of \textbf{$\bm{\kappa}$-invariant level $\bm{k}$
weights} of $G$:
\[
(\Lambda^{\ast})_{k}^{\kappa}:=(\Lambda_{+}^{\ast})^{\kappa}\cap\mathfrak{A}_{(\kappa),k}^{\ast}.
\]

\item The finite subgroup:
\[
T_{k+\mathsf{h}_{(\kappa)}^{\vee}}^{\kappa}:=\tilde{B}_{k+\mathsf{h}_{(\kappa)}^{\vee}}^{\sharp}\left((\Lambda^{\ast})^{\kappa}\right)/\Lambda_{(\kappa)}\subseteq T_{(\kappa)},
\]
where $\mathsf{h}_{(\kappa)}^{\vee}=1+\langle\rho,\tilde{B}^{\sharp}(\theta_{(\kappa),\mathrm{l}})\rangle$
is the dual Coxeter number of $\mathfrak{R}_{(\kappa)}$, and where
$T_{k+\mathsf{h}_{(\kappa)}^{\vee}}^{\kappa,\mathrm{reg}}/W^{\kappa}$
consists of the elements:
\[
s_{\lambda}:=\exp_{\mathfrak{g}_{(\kappa)}}\left(\tilde{B}_{k+\mathsf{h}_{(\kappa)}^{\vee}}^{\sharp}(\lambda+\rho)\right),\ \ \lambda\in(\Lambda^{\ast})_{k}^{\kappa}.
\]

\item The shifted action $\bullet_{k}$ of $W_{\mathrm{aff}}^{(\kappa)}=\Lambda_{(\kappa)}\rtimes W^{\kappa}$
on $(\mathfrak{t}^{\kappa})^{\ast}$, specified by the assignment:
\[
w\bullet_{k}\lambda=\begin{cases}
w\cdot(\lambda+\rho)-\rho, & w\in W^{\kappa};\\
\lambda-\tilde{B}_{k+\mathsf{h}_{(\kappa)}^{\vee}}^{\flat}(w), & w\in\Lambda_{(\kappa)},
\end{cases}
\]
for any $\lambda\in(\mathfrak{t}^{\kappa})^{\ast}$.

\end{itemize} Note that the shifted action of $W_{\mathrm{aff}}^{(\kappa)}$
is generated by reflections in the affine hyperplanes:
\[
\mathcal{H}_{\alpha,m}^{(\kappa,k)}=\left\{ \lambda\in(\mathfrak{t}^{\kappa})^{\ast}\mbox{ }\big|\mbox{ }\langle\tilde{\alpha},\tilde{B}_{k+\mathsf{h}_{(\kappa)}^{\vee}}^{\sharp}(\lambda+\rho)\rangle=m\right\} ,\ \ \tilde{\alpha}\in\mathfrak{R}_{(\kappa)},m\in\mathbb{Z},
\]
and that the elements $s_{\lambda}$ are in bijection with those of
\[
(\Lambda^{\ast})_{k+\mathsf{h}_{(\kappa)}^{\vee},\mathrm{reg}}^{\kappa}=(\Lambda^{\ast})_{k}^{\kappa}+\rho.
\]
An element $\lambda\in(\Lambda_{+}^{\ast})^{\kappa}$ has a trivial
stabilizer under $W_{\mathrm{aff}}^{(\kappa)}$ if and only if $\lambda\in W_{\mathrm{aff}}^{(\kappa)}\bullet_{k}(\Lambda^{\ast})_{k}^{\kappa}$,
and is otherwise lying on some shifted affine wall $\mathcal{H}_{\alpha,m}^{(\kappa,k)}$.
Lastly, using the Jantzen character formula, we have for any $\lambda\in\lambda\in(\Lambda_{+}^{\ast})^{\kappa}$
and $w\in W_{\mathrm{aff}}^{(\kappa)}$:
\[
\tilde{\chi}_{w\bullet_{k}\lambda}^{\kappa}=(-1)^{l(w)}\tilde{\chi}_{\lambda}^{\kappa}.
\]

\begin{defn}
With the notations above, define the \textbf{level $\bm{k}$ twining
fusion ring} as the quotient:
\[
\tilde{R}_{k}^{(\kappa)}(G):=\tilde{R}^{(\kappa)}(G)/\tilde{I}_{k}^{(\kappa)}(G),
\]
where $\tilde{I}_{k}^{(\kappa)}(G)$ is the level $k$ fusion ideal:
\[
\tilde{I}_{k}^{(\kappa)}(G):=\left\{ f\in\tilde{R}^{(\kappa)}(G)\mbox{ }\big|\mbox{ }f(t)=0\mbox{, for all }t\in\pi^{-1}(s_{\lambda}),\lambda\in(\Lambda^{\ast})_{k}^{\kappa}\right\} ,
\]
with $\pi:T^{\kappa}\rightarrow T_{(\kappa)}$ the projection.
\end{defn}

The $\mathbb{Z}$-linear maps $(\cdot)^{\ast}:\tilde{R}^{(\kappa)}(G)\rightarrow\tilde{R}^{(\kappa)}(G)$
and $\mathrm{Tr}_{0}:\tilde{R}^{(\kappa)}(G)\rightarrow\mathbb{Z}$
descend to the ring $\tilde{R}_{k}^{(\kappa)}(G)$, which has $\{\tilde{\chi}_{\mu}^{\kappa}\}_{\mu\in(\Lambda^{\ast})_{k}^{\kappa}}$
as a basis over $\mathbb{Z}$, and:
\[
\tilde{R}_{k}^{(\kappa)}(G)\simeq\mathbb{Z}\left[(\Lambda^{\ast})_{k}^{\kappa}\right].
\]
The isomorphism $\tilde{R}^{(\kappa)}(G)\simeq R\left(G_{(\kappa)}\right)$
induces an isomorphism of $\tilde{R}^{(\kappa)}(G)$-modules $\tilde{I}_{k}^{(\kappa)}(G)\simeq I_{k}\left(G_{(\kappa)}\right)$,
with $I_{k}\left(G_{(\kappa)}\right)\subseteq R\left(G_{(\kappa)}\right)$
the level $k$ fusion ideal of $G_{(\kappa)}$. We thus obtain a complete
characterization of $\tilde{R}_{k}^{(\kappa)}(G)$ for $k\in\mathbb{Z}_{>0}$:
\begin{prop}
Let $G$ be compact, connected, simply connected and simple, $\kappa\in\mathrm{Out}(G)$
a diagram automorphism and $G_{(\kappa)}$ the corresponding orbit
Lie group. The twining fusion ring $\tilde{R}_{k}^{(\kappa)}(G)$
at level $k\in\mathbb{Z}_{>0}$ satisfies the following properties:\begin{enumerate}

\item It is isomorphic to the fusion ring of $G_{(\kappa)}$ at level
$k$:
\[
\tilde{R}_{k}^{(\kappa)}(G)\simeq R_{k}\left(G_{(\kappa)}\right).
\]

\item Under the canonical projection $\Phi:\tilde{R}^{(\kappa)}(G)\rightarrow\tilde{R}_{k}^{(\kappa)}(G)$,
the images of basis elements $\{\tilde{\chi}_{\mu}^{\kappa}\}_{\mu\in(\Lambda^{\ast})_{+}^{\kappa}}$
of $\tilde{R}^{(\kappa)}(G)$ are given by: 
\[
\Phi\left(\tilde{\chi}_{\mu}^{\kappa}\right)=\begin{cases}
(-1)^{l(w)}\tilde{\chi}_{w\bullet_{k}\mu}^{\kappa}, & \mbox{ if }w\bullet_{k}\mu\in\Lambda_{k}^{\ast}\mbox{ for }w\in W_{\mathrm{aff}}^{(\kappa)}\smallsetminus\{1\};\\
0, & \mbox{ if }(W_{\mathrm{aff}}^{(\kappa)})_{\mu}\ne\{1\}.
\end{cases}
\]
The elements $\left\{ \Phi\left(\tilde{\chi}_{\mu}^{\kappa}\right)\right\} _{\mu\in(\Lambda^{\ast})_{k}^{\kappa}}$
constitute the spectrum of $\tilde{R}_{k}^{(\kappa)}(G)$.

\item For any $\lambda,\mu\in(\Lambda^{\ast})_{k}^{\kappa}$, the
fusion product on $\tilde{R}_{k}^{(\kappa)}(G)$ is given by:
\[
\tilde{\chi}_{\lambda}^{\kappa}\cdot\tilde{\chi}_{\mu}^{\kappa}=\sum_{\nu\in(\Lambda^{\ast})_{k}^{\kappa}}N_{\lambda\mu\nu^{\ast}}^{(\kappa,k)}\tilde{\chi}_{\nu}^{\kappa},
\]
where the fusion coefficients $N_{\lambda\mu\nu^{\ast}}^{(\kappa,k)}=\mathrm{Tr}_{0}(\tilde{\chi}_{\lambda}^{\kappa}\cdot\tilde{\chi}_{\mu}^{\kappa}\cdot\tilde{\chi}_{\nu^{\ast}}^{\kappa})$
are given by:
\[
N_{\lambda\mu\nu^{\ast}}^{(\kappa,k)}=\frac{1}{\Big|T_{k+\mathsf{h}_{(\kappa)}^{\vee}}^{\kappa}\Big|}\sum_{\alpha\in(\Lambda^{\ast})_{k}^{\kappa}}|\widetilde{J}_{0}(s_{\alpha})|^{2}\tilde{\chi}_{\lambda}^{\kappa}(s_{\alpha})\tilde{\chi}_{\mu}^{\kappa}(s_{\alpha})\tilde{\chi}_{\nu^{\ast}}^{\kappa}(s_{\alpha}).
\]

\end{enumerate}
\end{prop}

\noindent\emph{Outline of proof.} Statement (1) follows from the
discussion preceding the proposition. The same remarks are used to
modify the proofs of propositions 8.3 and 9.4 of \cite{[Be93]} to
obtain statement (2), and statement (3) is obtained by modifying the
proof of lemma 9.7 in \cite{[Be93]} (see also \cite[Prop. A.8]{[Mein12]}).
\qed

\hfill

To summarize, given a compact 1-connected simple Lie group $G$, the
representation and fusion rings associated to an automorphism $\kappa\in\mathrm{Out}(G)$
are isomorphic to those of the orbit Lie group $G_{(\kappa)}$. This
raises several questions, some of which we hope to address in future
work. For instance, it would be desirable to have more conceptual
proofs of the results in this section, and to see how they generalize
to the case of $G$ non-simply connected.

\bibliographystyle{amsplain}
\bibliography{Twisted_DH_Measure_Bibtex}

\end{document}